# Exchangeable pairs and Poisson approximation

**Sourav Chatterjee, Persi Diaconis and Elizabeth Meckes**

*Stanford University*
*e-mail:*
`souravc@stat.stanford.edu; diaconis@math.stanford.edu; meckes@math.stanford.edu`

**Abstract:** This is a survey paper on Poisson approximation using Stein's method of exchangeable pairs. We illustrate using Poisson-binomial trials and many variations on three classical problems of combinatorial probability: the matching problem, the coupon collector's problem, and the birthday problem. While many details are new, the results are closely related to a body of work developed by Andrew Barbour, Louis Chen, Richard Arratia, Lou Gordon, Larry Goldstein, and their collaborators. Some comparison with these other approaches is offered.

## 1. Introduction

Charles Stein has introduced a general method for proving limit theorems with explicit error terms, the essential idea of which is the notion of a characterizing operator. Given a probability measure $\mathbb{P}_o$ on a space $\Omega$ with expectation $\mathbb{E}_o$, a characterizing operator $T_o$ is an operator on a suitable function space on $\Omega$ with the properties that:

1. $\mathbb{E}_o T_o \equiv 0$.
2. If $\mathbb{E}$ is expectation for a different probability on $\Omega$ and $\mathbb{E}T_o \equiv 0$, then $\mathbb{E} = \mathbb{E}_o$.

The idea is then to prove that $\mathbb{E} \doteq \mathbb{E}_o$ by showing that $\mathbb{E}T_o \doteq 0$. To do this, Stein has introduced a method which he calls the method of exchangeable pairs.

In this survey, we specialize his general approach to Poisson approximation. Very roughly, the idea is to show that a random variable $W$ has an approximate Poisson distribution by studying how a small stochastic change affects the law of W. An appropriate such change is often easily found by constructing a reversible Markov chain resulting in an exchangeable pair $(W, W')$. Stein's general approach is developed for the Poisson setting in section 2. This development is somewhat technical; as motivation for the machinery, we offer the following rough overview. The aim is to prove that $W = \sum X_i$ has an approximate Poisson distribution, where the $X_i$ are (perhaps dependent) indicators. First, a probabilistic construction is used to create a second random variable $W'$ such that $(W, W')$ is exchangeable. This construction is often one step in a reversible Markov chain with $W$ as stationary distribution. Then the machinery from section 2 gives the following bound:

$$|\mathbb{P}(W \in A) - Poi_\lambda(A)| \leq (0.7)\lambda^{-1/2}(E_1 + E_2)$$

with

$$
\begin{aligned}
E_1 &= \mathbb{E}|\lambda - c\mathbb{P}(W' = W + 1|\{X_i\})| \\
E_2 &= \mathbb{E}|W - c\mathbb{P}(W' = W - 1|\{X_i\})|
\end{aligned}
$$

with $c$ a parameter chosen to minimize the error terms. Thus the error will be small provided

$$\mathbb{P}(W' = W + 1|\{X_i\}) \doteq \frac{\lambda}{c}$$

and

$$\mathbb{P}(W' = W - 1|\{X_i\}) \doteq \frac{W}{c}.$$

In applications, $W = \sum X_i$ with the $X_i$ indicators of rare events. The $X_i$ are thus mostly zero and $W$ of them are one. The example to keep in mind is when a point of the sequence is chosen at random, and that $X_i$ is changed to its opposite, giving $W'$. Heuristically, since most of the $X_i$ are 0, there are constants $a$ and $b$ such that

$$
\begin{aligned}
\mathbb{P}(W' = W + 1|\{X_i\}) &\doteq a \\
\mathbb{P}(W' = W - 1|\{X_i\}) &\doteq bW.
\end{aligned}
$$

By the symmetry of the exchangeable pair,

$$
\begin{aligned}
a &\doteq \mathbb{E}[\mathbb{P}(W' = W + 1|\{X_i\})] \\
&= \mathbb{P}(W' = W + 1) \\
&= \mathbb{P}(W = W' + 1) \\
&= \mathbb{E}[\mathbb{P}(W' = W - 1|\{X_i\})] \\
&\doteq \mathbb{E}(bW) \\
&= b\lambda.
\end{aligned}
$$

Thus choosing $c = \frac{\lambda}{a}$ gives $b = \frac{1}{c}$, and makes the error terms small. In many cases, just changing one $X_i$ to its opposite doesn't quite work to produce an exchangeable pair, but most of our exchangeable pairs are constructed in ways which are rather similar to this. The reader will see rigorous versions of these heuristics starting in section 3.

The contents of the rest of the paper are as follows. Section 2 develops a general bound. In section 3, we give a first example of using Stein's method: Poisson-binomial trials. The method is shown to give close to optimal results with explicit error bounds. In sections 4-6, the method is applied to three classical probability problems: the matching problem, the birthday problem, and the coupon collector's problem. Sections 7 and 8 give multivariate and Poisson process approximations.

The Poisson heuristic and the Chen-Stein approach to Poisson approximation have been actively developed since 1970. Substantial contributions from Louis Chen, Andrew Barbour and their many coauthors are surveyed in sections 9 and 10. These approaches share common roots with the present exchangeable pair approach, and all have different strengths and weaknesses; we attempt some comparison. Section 11 lists some interesting open problems.

**Acknowledgements.** This paper is based on notes from a seminar at Stanford University in the winter of 2004. We thank the other speakers and participants, particularly Susan Holmes and Charles Stein. We also thank Geir Helleloid for a careful reading of the manuscript.

## 2. The Method of Exchangeable Pairs and Poisson Approximation

Let $W$ be an integer-valued random variable defined on a probability space $(\Omega, \mathcal{A}, \mathbb{P})$. Let $\mathcal{X}$ be the bounded measurable functions on $\Omega$. We write $\mathbb{E}f = \int f(\omega)\mathbb{P}(d\omega)$ and think of $\mathbb{E}$ as a linear map, $\mathbb{E} : \mathcal{X} \to \mathbb{R}$. Throughout we suppose that $\mathbb{E}(W) = \lambda < \infty$. Let $\mathcal{X}_o$ be the bounded real-valued functions on $\mathbb{N} = \{0, 1, 2, \ldots\}$, and let $\mathbb{E}_o : \mathcal{X}_o \to \mathbb{R}$ be expectation with respect to $Poi_\lambda$ measure on $\mathbb{N}$. The random variable $W$ allows us to define a map $\beta : \mathcal{X}_o \to \mathcal{X}$ by

$$\beta f(\omega) = f(W(\omega)). \tag{1}$$

Heuristically, $W$ has an approximate $Poi_\lambda$ distribution if $\mathbb{E}_o(f) \doteq \mathbb{E}(f(W))$; i.e., $\mathbb{E}_o \doteq \mathbb{E}\beta$. This is equivalent to saying that the following diagram approximately commutes:

$$\tag{2}$$

Stein constructs a symmetric probability $\mathbb{Q}$ on $\Omega \times \Omega$ with margins $\mathbb{P}$ (i.e., $\mathbb{Q}(A, B) = \mathbb{Q}(B, A)$ and $\mathbb{Q}(A, \Omega) = \mathbb{P}(A)$), which gives an exchangeable pair. It is used to define the following enlarged diagram:

(3)

Stein's lemma, developed and proved in lemma (2) below, shows in a precise sense that if the left square approximately commutes, then the triangle approximately commutes. This leads to explicit bounds on Poisson approximation.

We begin by constructing the top row of the diagram. We define the characterizing or Stein operator $T_o$ for the $Poi_\lambda$ distribution by

$$T_o f(j) = \lambda f(j+1) - j f(j) \tag{4}$$

and we let $\mathcal{F}_o \subseteq \mathcal{X}_o$ be the functions $f : \mathbb{N} \to \mathbb{R}$ such that $T_o f$ is a bounded function. Note that $\mathcal{F}_o$ contains all of the functions $f$ on $\mathbb{N}$ such that $f(n) = 0$ eventually, thus it is rich enough for our purposes.

A simple calculation verifies that $\mathbb{E}_o T_o = 0$; i.e., if $X \sim Poi_\lambda$, then

$$\mathbb{E}(\lambda f(X+1) - X f(X)) = 0$$

for every bounded function $f$.

**Remark:** $T_o$ is called characterizing because it also has the property that if $p$ is a probability on $\mathbb{N}$ with the property that $\sum_{j=0}^{\infty} (\lambda f(j+1) - j f(j)) p(j) = 0$ for every bounded function $f : \mathbb{N} \to \mathbb{R}$, then $p$ is the Poisson distribution with parameter $\lambda$. To see this, let $f = \delta_k$. This yields the equation $\lambda p(k-1) = k p(k)$, which leads to a recursion relation that describes $Poi_\lambda$.

Next we will define a map $U_o$, which is almost an inverse to $T_o$. Define:

$$U_o f(j) = \frac{(j-1)!}{\lambda^j} \sum_{k=0}^{j-1} \frac{\lambda^k}{k!} (f(k) - \mathbb{E}_o f). \tag{5}$$

It is easy to check that

$$T_o U_o f(j) = f(j) - \mathbb{E}_o f. \tag{6}$$

Thus $U_o$ is inverse to $T_o$ on $ker(\mathbb{E}_o)$. The following lemma, proved by Barbour and Eagleson in [7], gives bounds on expressions involving $U_o$.

**Lemma 1.** *1. For $f \in \mathfrak{X}_o$ with $0 \leq f \leq 1$ and $j \in \mathbb{N}$, the map $U_o$ of (5) satisfies*

$$
\begin{aligned}
|U_o f(j)| &\leq \min(1, 1.4\lambda^{-1/2}) \\
|U_o f(j+1) - U_o f(j)| &\leq \frac{1 - e^{-\lambda}}{\lambda}
\end{aligned}
$$

*2. Let $\delta_o(j) = \begin{cases} 1 & j = 0, \\ 0 & otherwise. \end{cases}$*

*For $f(j) = \delta_o(j) - e^{-\lambda}$,*

$$
|U_o f(j)| \leq \frac{1 - e^{-\lambda}}{\lambda}.
$$

A reviewer has pointed to a sharper bound in [38]. We have not seen this and will use lemma (1) in what follows. This lemma and equation (6) show that $U_o f \in \mathcal{F}_o$ for every $f \in \mathfrak{X}_o$. To complete the top row, let $i$ denote the map $\mathbb{R} \to \mathfrak{X}_o$ which associates to each constant $c$ the map on $\mathbb{N}$ which takes the value $c$ for each $n$.

To construct the bottom row of diagram (3), define

$$
\mathcal{F} = \{ f : \Omega \times \Omega \to \mathbb{R} : f \text{ bounded, measurable, and} f(\omega, \omega') = -f(\omega', \omega) \}.
$$

Using a probability $\mathbb{Q}$ on $\Omega \times \Omega$ as discussed above, define the operator $T$ by:

$$
Tf(\omega) = \mathbb{E}_{\mathbb{Q}}(f(\omega, \omega') | \omega). \tag{7}
$$

Observe that $Tf$ is a bounded function on $\Omega$ for $f \in \mathcal{F}$. Further,

$$
\begin{aligned}
\mathbb{E}Tf &= \int f(\omega, \omega') \mathbb{Q}(d\omega, d\omega') \\
&= \int f(\omega', \omega) \mathbb{Q}(d\omega, d\omega') \\
&= -\int f(\omega, \omega') \mathbb{Q}(d\omega, d\omega')
\end{aligned}
$$

where the first equality is by the symmetry of $\mathbb{Q}$ and the second is by the anti-symmetry of $f$. Thus $\mathbb{E}T = 0$.

Finally, define $\alpha : \mathcal{F}_o \to \mathcal{F}$ to be any linear map, for example, $\alpha f(\omega, \omega') = f(W(\omega)) - f(W(\omega'))$. Stein's lemma is true regardless of what $\alpha$ is, so we choose $\alpha$ to work well with the problem at hand. In applications, $\alpha$ is often a localized version of the example given above; see proposition (3) below for an example.

We can now state and prove:

**Lemma 2 (Stein).** *Suppose that in the following diagram of linear spaces and linear maps, $\mathbb{E}T = 0$ and $T_o U_o = Id - i\mathbb{E}_o$.*

$$(8)$$

*Then*

$$\mathbb{E}\beta - \mathbb{E}_o = \mathbb{E}(\beta T_o - T\alpha)U_o \qquad (9)$$

*Proof.* We have:

$$
\begin{aligned}
0 &= \mathbb{E}T\alpha U_o \\
&= \mathbb{E}(T\alpha - \beta T_o)U_o + \mathbb{E}\beta T_o U_o \\
&= \mathbb{E}(T\alpha - \beta T_o)U_o + \mathbb{E}\beta(Id - i\mathbb{E}_o) \\
&= \mathbb{E}(T\alpha - \beta T_o)U_o + \mathbb{E}\beta - \mathbb{E}_o
\end{aligned}
$$

as desired. $\square$

**Remarks:**

1. Stein's lemma makes precise the sense in which the Poisson approximation $\mathbb{E}\beta \doteq \mathbb{E}_o$ holds provided $\beta T_o - T\alpha$ is small; this expression must have small expectation on the range of $U_o$. A good choice of $\alpha$ can put this expression into a useful form (see proposition (3) below).

2. In what follows, the equality in Stein's lemma is used to bound the error term $\mathbb{E}(T\alpha - \beta T_o)U_o$ with analysis. The same equality can also be used to simulate the error term using the exchangeable pair as a reversible Markov chain. This approach is developed in Stein et al. (2003).

3. Stein's method is often used to study a random variable $W$ which is a sum of indicators $X_i$. In the previous development, the operator $T$ was defined as conditional on $W$. The lemma holds if $T$ is defined by conditioning on any $\sigma$-algebra with respect to which $W$ is measurable; in particular, one can condition on $\{X_i\}$. An application of Jensen's inequality shows that conditioning on $W$ alone rather than conditioning on $\{X_i\}$ always gives better results. However, conditioning on only $W$ can make bounding the error terms more difficult, and in practice we often condition on all of the $X_i$.

4. One reason for presenting the diagrammatic view of the subject is that it applies, in essentially the same way, to any type of approximation. The fact that we're approximating by a Poisson random variable enters into the diagram only in the top row, as the characterizing operator is the one corresponding to the $Poi_\lambda$ distribution. To use Stein's method in other contexts, one needs only to find the appropriate characterizing operator; the rest of the picture remains unchanged. Holmes (in chap. 3 of [23]) gives bounds for general discrete measures.

5. Note that for lemma (2) and the bounds to follow, it is only required that $W$ and $W'$ are exchangeable. While this is most easily achieved by making $\mathbb{Q}$ symmetric, examples of Fulman [28] and Rinott-Rotar [45] show that other constructions are possible.

The following proposition gives us a more workable form for the error term.

**Proposition 3.** *Let $W$ be a random variable on $(\Omega, \mathcal{A}, \mathbb{P})$ and $\mathbb{Q}$ a symmetric probability on $\Omega \times \Omega$ with margins $\mathbb{P}$. Let $W = W(\omega)$ and $W' = W(\omega')$. Then*

1. *$\beta T_o f = \lambda f(W + 1) - W f(W)$ for $f \in \mathfrak{X}_o$.*
2. *Let $\alpha f(\omega, \omega') = c f(W') \mathbb{I}_{\{W'-1\}}(W) - c f(W) \mathbb{I}_{\{W-1\}}(W')$ for a constant $c$. Then for $f \in \mathcal{F}_o$,*

$$T\alpha f(\omega) = c f(W+1)\mathbb{Q}(W' = W+1 \,|\, \omega) - c f(W)\mathbb{Q}(W' = W-1 \,|\, \omega). \quad (10)$$

**Remarks:** For this choice of $\alpha$, the error term from Stein's lemma will be small if the exchangeable pair can be chosen so that

$$\mathbb{Q}(W' = W - 1 | \omega) \quad \approx \quad \frac{W}{c} \qquad\qquad (11)$$

$$\mathbb{Q}(W' = W + 1 | \omega) \quad \approx \quad \frac{\lambda}{c} \qquad\qquad (12)$$

for some $c$. One then uses this $c$ in defining $\alpha$ to obtain cancellation in the error term, and then uses the bounds in lemma (1) to bound the error term. The examples which follow show that (11) and (12) often hold for natural choices of $\mathbb{Q}$, for some $c$. We observe that if (11) and (12) hold, the ratio method proposed by Stein [52] is also a way to approach Poisson approximation.

## 3. Poisson-Binomial Trials

Let $X_1, X_2, ..., X_n$ be independent $\{0, 1\}$-valued random variables with $\mathbb{P}(X_i = 1) = p_i$ and $\mathbb{P}(X_i = 0) = 1 - p_i$, such that $\sum_{i=1}^{n} p_i = \lambda$. To put this example into the framework of section (2), let $\Omega = \{0, 1\}^n$, with $\mathbb{P}(\omega_1, ..., \omega_n) = \prod_{i=1}^{n} p_i^{\omega_i}(1 - p_i)^{1-\omega_i}$, and let $\mathfrak{X}$ be the bounded functions on $\Omega$. Let $W(\omega) = \sum_{i=1}^{n} \omega_i$. One way to build a symmetric probability $\mathbb{Q}$ on $\Omega \times \Omega$ is the following probabilistic construction: choose an index $I$ uniformly in $\{1,...,n\}$, and let $\epsilon_I = 1$ with probability $p_I$ and $\epsilon_I = 0$ with probability $1 - p_I$. Given $\omega \in \Omega$, set $\omega'_i = \begin{cases} \omega_i & \text{if } i \neq I, \\ \epsilon_I & \text{if } i = I. \end{cases}$ This constructs a probability $\mathbb{Q}$ on $\Omega \times \Omega$: assign to the pair $(\omega, \omega')$ the probability of choosing $\omega$ from $\Omega$ according to $\mathbb{P}$ and then going from $\omega$ to $\omega'$ by the process described. It is clear from the construction that $\mathbb{Q}$ is symmetric and has margins $\mathbb{P}$. From these definitions,

$$
\begin{aligned}
\mathbb{Q}(W' = W - 1 | \omega) &= \frac{1}{n} \sum_{i=1}^{n} \omega_i (1 - p_i) \\
&= \frac{W - \sum_{i=1}^{n} p_i \omega_i}{n}.
\end{aligned}
$$

$$
\begin{aligned}
\mathbb{Q}(W' = W + 1 | \omega) &= \frac{1}{n} \sum_{i=1}^{n} (1 - \omega_i) p_i \\
&= \frac{\lambda - \sum_{i=1}^{n} p_i \omega_i}{n}.
\end{aligned}
$$

Combining these calculations with (10) and choosing $c = n$ gives

$$
\begin{aligned}
(\beta T_o - T\alpha) f(\omega) &= [\lambda f(W+1) - W f(W)] - [f(W+1)(\lambda - \sum_{i=1}^{n} p_i \omega_i) - f(W)(W - \sum_{i=1}^{n} p_i \omega_i)] \\
&= (f(W+1) - f(W)) \sum_{i=1}^{n} p_i \omega_i.
\end{aligned}
$$

As we consider $f = U_o g$ for functions $g$ with $0 \leq g \leq 1$, lemma (1) yields the following:

**Theorem 4.** *Let $X_1, X_2, ..., X_n$ be independent $\{0,1\}$-valued random variables with $\mathbb{P}(X_i = 1) = p_i$ and $\sum_{i-1}^{n} p_i = \lambda$. Then, if $W = \sum_{i=1}^{n} X_i$,*

$$
\|\mathcal{L}(W) - Poi_\lambda\|_{TV} \leq \left( \frac{1 - e^{-\lambda}}{2\lambda} \right) \sum_{i=1}^{n} p_i^2.
$$

**Remark:** The factor of 2 in the denominator arises because

$$
\|\mu - \nu\|_{TV} = \sup_A |\mu(A) - \nu(A)| = \frac{1}{2} \sup_f \left| \int_A d\mu - \int_A d\nu \right|,
$$

where the first supremum is taken over all measurable sets $A$ and the second is over all measurable funtions $f$ such that $\|f\|_\infty \leq 1$.

**Example 1:** Let $p_i = \frac{\lambda}{n}$. For $\lambda$ fixed and $n > \lambda$, the bound becomes $\frac{(1-e^{-\lambda})\lambda}{2n}$. This bound was derived by a different argument by Barbour and Hall. A comparison with other available bounds and approximations for the i.i.d. case is in Kennedy and Quine [35], where there are also extensive references given. The bound that we give is sharp for small $\lambda$ and can be improved for large $\lambda$.

**Example 2:** Let $p_i = \frac{1}{i}$ for $1 \le i \le n$. Thus $\lambda = \log n - \gamma + O\left(\frac{1}{n}\right)$. The bound becomes $\left(\dfrac{(1 - \frac{e^\gamma}{n})\pi^2}{12 \log n}\right)\left(1 + O\left(\dfrac{1}{n}\right)\right)$. Here, the factor of $\lambda$ in the denominator saves the day.

**Remarks:** Under the conditions of theorem (4), there are methods for computing the exact distribution of $W$ and a host of approximations that differ from the Poisson. See Stein [49] or Percus and Percus [42] for further discussion. An extensive collection of random variables which can be represented as sums of independent binary indicators is in Pitman [43].

## 4. The Matching Problem

Because of its appearance in Montmort [18], the matching problem is one of the oldest problems in probability. It asks for the distribution of the number of fixed points in a random permutation. Takács [53] gives an extensive history. To put this problem into our set up, let $S_n$ be all $n!$ permutations of $n$ objects with $\mathbb{P}(\sigma) = \frac{1}{n!}$, and let $W(\sigma) = |\{i : \sigma(i) = i\}|$; that is, $W(\sigma)$ is the number of fixed points of $\sigma$. Note that $W = \sum_{i=1}^{n} X_i$ where $X_i$ indicates whether $\sigma$ fixes $i$.

This approach of writing a random variable as a sum of $\{0, 1\}$-valued random variables, called the method of indicators, is one of our main tools for computing expectations of $\mathbb{N}$-valued random variables. Using this representation of $W$, it is easy to see that $\mathbb{E}W = 1$. Build an exchangeable pair $(\sigma, \sigma')$ by choosing $\sigma \in S_n$ according to $\mathbb{P}$ and then choosing $\sigma'$ given $\sigma$ by following $\sigma$ with a random transposition $\tau$ chosen uniformly among the non-trivial transpositions in $S_n$. Then

$$\mathbb{Q}(W' = W - 1 \,|\, \sigma) \;=\; \frac{2W(\sigma)(n - W(\sigma))}{n(n-1)}$$

$$\mathbb{Q}(W' = W + 1 \,|\, \sigma) \;=\; \frac{2(n - W(\sigma) - 2a_2(\sigma))}{n(n-1)}$$

where $a_2(\sigma)$ is the number of transpositions which occur when $\sigma$ is written as a product of disjoint cycles. The first expression is calculated by observing that if $W' = W - 1$, then $\tau$ must choose one of the fixed points and switch it with something which is not a fixed point. The second calculation is done by similar considerations.

Define $\alpha$ to be the localized operator as in proposition (3) and choose $c = \frac{n-1}{2}$.

Then, for $g : \mathbb{N} \to [0,1]$,

$$
\begin{aligned}
(\beta T_o - T\alpha)U_o g(\sigma) =& [U_o g(W+1) - W \cdot U_o g(W)] \\
&- \left[ U_o g(W+1) \frac{(n-W-2a_2)}{n} - U_o g(W) \frac{W(n-W)}{n} \right] \\
=& U_o g(W+1) \frac{(W+2a_2)}{n} - U_o g(W) \frac{W^2}{n}
\end{aligned}
\tag{13}
$$

Using the method of indicators, it is easy to check that $\mathbb{E}(W^2) = 2$ and $\mathbb{E}(2a_2) = 1$. Further, for $\|g\|_\infty \le 1$, $\|U_o g\|_\infty \le 1$ by lemma (1). Putting all of this together proves

**Theorem 5.** *Let $W$ be the number of fixed points in a randomly chosen permutation on $n$ letters. Then*

$$
\|\mathcal{L}(W) - Poi_1\|_{TV} \le \frac{2}{n}.
$$

**Remarks:**

1. In this example, the bound is not sharp. As is well known, $\|\mathcal{L}(W) - Poi_1\|_{TV} \le \frac{2^n}{n!}$. It is interesting to speculate on just where sharp control is lost. We have proceeded by equality through (13).

2. The super exponential bound of remark (1) is an algebraic accident. This can be seen by considering the number of fixed points in only part of the permutation, that is, considering $W_k(\sigma) = |\{1 \le i \le k : \sigma(i) = i\}|$ for some $k < n$. Essentially the same argument shows that $W_k \stackrel{.}{\sim} Poi_{\frac{k}{n}}$. In section 8, we show that the fixed point process (scaled to [0,1]) converges to a Poisson process of rate 1. Here, the error may be seen to be of order $\frac{1}{n}$.

*Generalized matching*

Consider now the more general problem of fixed points of permutations of the set $A = \{1, 1, ..., 1, 2, ..., 2, ..., k, ..., k\}$ where the number $i$ appears $l_i$ times. For example, if two ordinary decks of cards are shuffled, placed on a table and turned up simultaneously, one card at a time, a match is counted if two cards of the same number appear at the same time, without regard to suits. This is the matching problem with $k = 13$ and $l_i = 4$ for each $i$.

Let $|A| = n$ and let $S_n$ be the set of all permutations of the elements of $A$ with $\mathbb{P}(\sigma) = \frac{1}{n!}$ for each $\sigma$. Let $W(\sigma)$ be the number of fixed points of $\sigma$ and $W_i(\sigma)$ the number of times that $i$ is a fixed point of $\sigma$, thus $W = \sum_1^k W_i$. Let $W_{ij}(\sigma) = |\{i : \sigma(i) = j\}|$, that is, $W_{ij}$ is the number of times $\sigma$ sends $i$ to $j$. Finally, let

$$X_m^i(j) = \begin{cases} 1 & \text{if } \sigma(i_m) = j, \\ 0 & \text{otherwise,} \end{cases}$$

where $i_m$ is the $m$th $i$ in the set $A$.

Build an exchangeable pair in this situation in the same way as in the previous situation: follow a permutation $\sigma$ by a random transposition of the $n$ set elements. Then

$$\begin{aligned} \mathbb{Q}(W' = W + 1|\sigma) &= \frac{1}{\binom{n}{2}} \sum_{i=1}^{k} \sum_{j \neq i} W_{ji}(l_i - W_i - W_{ij}) \\ &= \frac{2}{n(n-1)} \sum_{i=1}^{k} \left[ l_i^2 - 2l_i W_i + W_i^2 - \sum_{i \neq j} W_{ji} W_{ij} \right], \end{aligned}$$

$$\begin{aligned} \mathbb{Q}(W' = W - 1|\sigma) &= \frac{1}{\binom{n}{2}} \sum_{i} W_i(n - W - l_i + W_i) \\ &= \frac{2}{n(n-1)} \left[ nW - \sum_{i=1}^{k} l_i W_i - W^2 + \sum_{i=1}^{k} W_i^2 \right] \end{aligned}$$

To calculate the first expression, make use of the fact that for fixed $i$, $\sum_{j \neq i} W_{ji} = l_i - W_i$. The first line follows as one determines the probability in question by, for each fixed $i$ and each $j \neq i$, counting the number of ways to choose a symbol $i$ with $\sigma(j) = i$ and a symbol $k \neq i$, $k \neq j$ with $\sigma(i) = k$. These will then be switched by a transposition, and the number of fixed points will have been increased by 1. The second calculation is done similarly.

In order to bound the error term, the following moments are needed:

$$\lambda = \mathbb{E}(W) = \mathbb{E}\left( \sum_{i=1}^{k} \sum_{m=1}^{l_i} X_m^i(i) \right) = \frac{1}{n} \sum_{i=1}^{k} l_i^2.$$

$$\mathbb{E}(W_i) = \mathbb{E}\left( \sum_{m=1}^{l_i} X_m^i(i) \right) = \frac{l_i^2}{n}.$$

$$\mathbb{E}(W_i^2) = \mathbb{E}\left( \left( \sum_{m=1}^{l_i} X_m^i(i) \right)^2 \right) = \frac{l_i^2[n + l_i^2 - 2l_i]}{n(n-1)}.$$

$$\mathbb{E}(W_{ij}W_{ji}) = \mathbb{E}\left(\sum_{p=1}^{l_i}\sum_{m=1}^{l_j} X_p^i(j)X_m^j(i)\right) = \frac{l_i^2 l_j^2}{n(n-1)} \quad (\text{for } i \neq j).$$

$$
\begin{aligned}
\mathbb{E}\left(W^2 - \sum_{i=1}^{k} W_i^2\right) &= \mathbb{E}\left(\sum_{i=1}^{k}\sum_{j\neq i} W_i W_j\right) \\
&= \mathbb{E}\left(\sum_{i=1}^{k}\sum_{j\neq i}\sum_{p=1}^{l_i}\sum_{m=1}^{l_j} X_p^i(i)X_m^j(j)\right) \\
&= \frac{1}{n(n-1)}\sum_{i=1}^{k}\sum_{j\neq i} l_i^2 l_j^2.
\end{aligned}
$$

Now choose (in analogy with the previous case) $c = \frac{n-1}{2}$ and make use of proposition (1) to estimate:

$$
\begin{aligned}
|\mathbb{E}f(W+1)[\lambda \;-\; & c\mathbb{Q}(W' = W+1|\sigma)]| \\
&\leq (1.4)\lambda^{-1/2}\frac{1}{n}\mathbb{E}\left|\sum_{i=1}^{k}\left[2l_iW_i - W_i^2 + \sum_{j\neq i} W_{ij}W_{ji}\right]\right| \\
&\leq \frac{1.4}{n\lambda^{1/2}}\left[\mathbb{E}\left(\sum_{i=1}^{k} l_iW_i\right) + \mathbb{E}\left(\sum_{i=1}^{k} W_i^2\right) + \mathbb{E}\left(\sum_{i=1}^{k}\sum_{j\neq i} W_{ji}W_{ij}\right)\right] \\
&= \frac{1.4}{n^2(n-1)\lambda^{1/2}}\sum_{i=1}^{k}\left((2n-4)l_i^3 + l_i^2\sum_{j=1}^{k} l_j^2\right) \\
&= \frac{1.4}{n^2(n-1)\lambda^{1/2}}\sum_{i=1}^{k}[(2n-4)l_i^3 + l_i^2(n\lambda)] \\
&\leq (1.4)\left[\left(\frac{\lambda^{3/2}}{n-1}\right) + \left(\frac{2}{n^2\lambda^{1/2}}\right)\sum_{i=1}^{k} l_i^3\right]
\end{aligned}
$$

for $f = U_o g$ with $g : \mathbb{N} \to [0,1]$.

For the other half of the error term,

$$
\begin{aligned}
|\mathbb{E}f(W)(c\mathbb{Q}(W' = W - 1|\sigma) - W)| &\leq \frac{1.4}{n\lambda^{1/2}} \left[ \mathbb{E}\left(\sum_{i=1}^{k} l_i W_i\right) + \mathbb{E}\left(W^2 - \sum_{i=1}^{k} W_i^2\right) \right] \\
&= \frac{1.4}{n^2(n-1)\lambda^{1/2}} \sum_{i=1}^{k} \left[ l_i^3(n-1) + \sum_{j\neq i} l_i^2 l_j^2 \right] \\
&\leq (1.4)\left[ \left(\frac{\lambda^{3/2}}{n-1}\right) + \left(\frac{1}{n^2\lambda^{1/2}}\right) \sum_{i=1}^{k} l_i^3 \right].
\end{aligned}
$$

for $f = U_o g$ and $g : \mathbb{N} \to [0, 1]$.

Putting these two estimates together with proposition (3) proves

**Theorem 6.** *Let $A = \{1, ..., 1, 2, ..., 2, ..., k, ..., k\}$, with the number $i$ appearing $l_i$ times, $|A| = n$. Let $W$ be the number of fixed points of a random permutation of the elements of $A$, and let $\mu = \sum_{i=1}^{k} l_i^3$. Then $\lambda := \mathbb{E}(W) = \frac{1}{n}\sum_{i=1}^{k} l_i^2$ and*

$$
\|\mathcal{L}(W) - Poi_\lambda\|_{TV} \leq (1.4)\left[ \frac{\lambda^{3/2}}{n-1} + \frac{3\mu}{2n^2\lambda^{1/2}} \right].
$$

**Remarks:**

1. In particular, if $l_i = l$ for each $i$, then the theorem gives:

$$
\|\mathcal{L}(W) - Poi_l\|_{TV} \leq \frac{(3.5)l^{3/2}}{n-1}.
$$

2. A celebrated paper of Kaplansky [33] first showed that $W$ had an approximate Poisson distribution under suitable restrictions on the $l_i$.

3. Stein's method may be used to prove Poisson approximation in the matching problem for some non-uniform distributions $\mu$ on $S_n$. In outline, the technique is similar: construct $\mathbb{Q}$ on $S_n \times S_n$ by choosing $\sigma$ according to $\mu$ and then choosing $\sigma'$ by making a random transposition and using the Metropolis algorithm to ensure that $\sigma'$ is also $\mu$-distributed.

4. Essentially the same techniques can be used to study approximate matches. Here, $W = |\{i : |i - \sigma(i)| \leq c\}|$ for some $c$.

5. Generalizing the problem in remark 4, one may study $\sum A_{i\pi(i)}$ where $A_{ij}$ is an $n \times n$ matrix with entries in $\{0, 1\}$. Similarly, one may study the number of permutations which violate a given set of restrictions at most $k$ times with this method, if the restrictions are given by a $\{0, 1\}$-matrix (that is, $\pi(i) = j$ is allowed if $A_{ij} = 1$). Classically, these problems are solved using rook theory or the extensive development connected to the permanent (see [21] for extensive references).

## 5. The Birthday Problem

For the classical birthday problem, $k$ balls are dropped independently into $n$ boxes, such that any ball is equally likely to fall into any of the $n$ boxes. How many balls should there be to have some specified chance of having a box with more than one ball in it? Let us set up the problem in the framework of Stein's method. Let $[n] = \{1, \ldots, n\}$, and let $\Omega = [n]^k$. An element $\omega = (b_1, \ldots, b_k) \in \Omega$ will represent an arrangement of the $k$ balls in the $n$ boxes: $b_i$ is the number of the box into which we put ball $i$. Thus the $b_i$ are i.i.d. and uniformly chosen from $[n]$. Denote this probability measure on $\Omega$ by $\mathbb{P}$. Let $M_m$ be the number of boxes containing exactly $m$ balls, and let $M_{m+}$ be the number of boxes containing $m$ or more balls. We will show that if $k = \theta\sqrt{n}$, then $W = M_{2+}$ is approximately Poisson distributed with parameter $\frac{\theta^2}{2}$.

Consider the random variables $M_2, W$ and $V = \sum_{i<j} \mathbb{I}\{b_i = b_j\}$. Clearly, $M_2 \leq W \leq V$. Also, they are all equal when $M_{3+} = 0$. Now $M_{3+} \leq \sum_{i<j<l} \mathbb{I}\{b_i = b_j = b_l\}$, since any box which contains more than two balls must contain one such triplet as in the sum. Thus,

$$\mathbb{E}(M_{3+}) \leq \sum_{i<j<l} \mathbb{P}(b_i = b_j = b_l) = \binom{k}{3}\frac{1}{n^2} \leq \frac{\theta^3}{6\sqrt{n}}.$$

Thus, $\mathbb{P}(M_{3+} \neq 0) = \mathbb{P}(M_{3+} \geq 1) \leq \frac{\theta^3}{6\sqrt{n}}$. In particular, $d_{TV}(\mathcal{L}_W, \mathcal{L}_V) \leq \mathbb{P}(W \neq V) \leq \frac{\theta^3}{6\sqrt{n}}$. Note that $\mathbb{E}V = \binom{k}{2}\frac{1}{n} \approx \frac{\theta^2}{2}$.

Build an exchangeable pair $(\omega, \omega')$ as follows: choose $\omega \in \Omega$ according to $\mathbb{P}$ and then choose an index $I \in [k]$ uniformly. Let $b_I^*$ be i.i.d. with the $\{b_i\}$, and let $\omega' = (b_1, \ldots, b_I^*, \ldots, b_k)$. Thus the exchangeable pair is formed by first distributing the $k$ balls independently and uniformly in the $n$ boxes, and then choosing a ball at random and choosing a new box for it at random. As before, let $W = W(\omega)$, $W' = W(\omega')$, and compute:

$$
\begin{aligned}
\mathbb{Q}(W' = W - 1 \,|\, \omega) &= \frac{2M_2(n - M_1 - 1)}{kn} \\
\mathbb{Q}(W' = W + 1 \,|\, \omega) &= \frac{M_1}{k}\left(\frac{M_1}{n} - \frac{1}{n}\right) + \frac{k - M_1 - 2M_2}{k}\left(\frac{M_1}{n}\right) \\
&= \frac{M_1(k - 2M_2 - 1)}{kn}.
\end{aligned}
$$

The first equality holds as in order to make $W$ go down by 1, we must choose a ball from a box that contains exactly two balls and put it into any box except either the one it came from or any of the boxes which started with exactly one ball. The second computation is similar.

By proposition (3) and choosing $c = \frac{k}{2}$, the error term for Poisson approximation is given by:

$$\frac{1}{2}\sup\{\mathbb{E}(Af(W+1) - Bf(W)) : f : \mathbb{N} \to [0,1]\}$$

where

$$A = \frac{\theta^2}{2} - \frac{M_1(k - 2M_2 - 1)}{2n}$$

and

$$B = W - \frac{M_2(n - M_1 - 1)}{n}.$$

Observe that $\theta = kn^{-1/2}$ implies that

$$A = \frac{k(k - M_1) + 2M_1M_2 + M_1}{2n} \geq 0.$$

Now $\mathbb{E}(k - M_1) = \mathbb{E}(W) \leq \mathbb{E}(V) \leq \theta^2/2$, and $\mathbb{E}(M_1M_2) \leq k\mathbb{E}(M_2) \leq k\theta^2/2$. Thus,

$$\mathbb{E}|A| \leq \frac{k}{2n}\left(\frac{\theta^2}{2} + \theta^2 + 1\right) = \left(\left(\frac{3}{4}\right)\theta^3 + \left(\frac{1}{2}\right)\theta\right)\left(\frac{1}{\sqrt{n}}\right).$$

Finally,

$$\mathbb{E}|B| = \mathbb{E}\left|W - M_2 + \frac{M_2(M_1 + 1)}{n}\right| = \mathbb{E}(M_{3+}) + \frac{(k+1)\theta^2}{2n} \leq \frac{2\theta^3}{3\sqrt{n}} + \frac{\theta^2}{2n}$$

since $\mathbb{E}(M_{3+}) \leq \frac{\theta^3}{6\sqrt{n}}$, as shown before. Applying lemma (1), proves

**Theorem 7.** *If we drop $k = \theta\sqrt{n}$ balls independently and uniformly into $n$ boxes, and let $W$ be the number of boxes containing at least two balls, then*

$$\|\mathcal{L}(W) - Poi_{\theta^2/2}\|_{TV} \leq \min\{1, \sqrt{2}\theta^{-1}\}\left[\frac{19\theta^3 + 6\theta}{12\sqrt{n}} + \frac{\theta^2}{2n}\right].$$

### *Triple matches*

We next consider a variation of the birthday problem: triple matches. As before, let $\Omega$ be the space of possible arrangements of $k$ balls in $n$ boxes with probability $\mathbb{P}$ given by independent and uniform placement of the balls. This time, the random variable $W(\omega)$ will denote the number of triple matches in $\omega$; i.e., $|\{(i,j,l) : 1 \leq i < j < l \leq k, b_i = b_j = b_l\}|$, where $b_i$ is the box in which we put ball $i$. One can write $W = \sum_{i<j<l} \mathbb{I}_{\{b_i = b_j = b_l\}}$, and from this representation, we get $\mathbb{E}W = \binom{k}{3}\frac{1}{n^2}$.

Let $M_i(\omega)$ be the number of boxes containing exactly $i$ balls, thus $W = \sum_{l=3}^{k} \binom{l}{3} M_l$. Let $X_i^l$ be the indicator that box $i$ contains exactly $l$ balls. Then

$$
\begin{aligned}
\mathbb{E}(M_l) &= \sum_{i=1}^{n} \mathbb{E}(X_i^l) \\
&= n \binom{k}{l} \left(\frac{1}{n}\right)^l \left(1 - \frac{1}{n}\right)^{k-l}.
\end{aligned}
$$

Construct an exchangeable pair in this case exactly as in the previous case: choose a ball at random and then choose a new box at random to put it in. We have:

$$
\begin{aligned}
\mathbb{Q}(W' = W + 1 | \omega) &= \frac{1}{kn}[M_1 M_2 + 2M_2(M_2 - 1)] \\
&= \frac{1}{kn}[M_1 M_2 + 2M_2^2 - 2M_2].
\end{aligned}
$$

This follows as the only ways to make the number of triples sharing a box go up by exactly one are either to choose a ball from a box with only one ball and move it to any of the boxes with two balls, or to choose a ball from a box with two balls and put it into one of the other boxes with two balls. No further possibilities exist, as removing a ball from a box with $m$ balls and putting it into a box with $l$ balls will increase $W$ by exactly one if and only if $\binom{m}{3} - \binom{m-1}{3} = \binom{l}{2} = \binom{l}{2} - 1$, and there are no solutions of this if $m > 3$. Similarly:

$$
\mathbb{Q}(W' = W - 1 | \omega) = \frac{1}{kn}[3M_3 M_0 + 3M_3 M_1].
$$

Using lemma (1) and proposition (3) together with the triangle inequality, we will need to estimate:

$$
\mathbb{E} \left| \frac{k(k-1)(k-2)}{6n^2} - \frac{c}{kn}[M_1 M_2 + 2M_2^2 - 2M_2] \right| \tag{14}
$$

$$
\mathbb{E} \left| \frac{c}{kn}[3M_0 M_3 + 3M_1 M_3] - W \right| \tag{15}
$$

for some choice of the parameter $c$.

Consider (14) first. Define new random variables $\widetilde{M_i} = M_i - \mu_i$ where $\mu_i = \mathbb{E} M_i$ and write everything in terms of these $\widetilde{M_i}$. We have:

$$
M_1 M_2 + 2M_2^2 - 2M_2 = \widetilde{M_1}\widetilde{M_2} + 2\widetilde{M_2}^2 + (\mu_1 + 4\mu_2 - 2)\widetilde{M_2} + \mu_2 \widetilde{M_1} + \mu_1 \mu_2 - 2\mu_2 + 2\mu_2^2.
$$

The deterministic part of this is

$$
\mu_1 \mu_2 - 2\mu_2 + 2\mu_2^2 = \frac{k(k-1)(n-1)^{k-2}[k(n-1)^{k-1} - 2n^{k-1} + k(k-1)(n-1)^{k-2}]}{2n^{2k-2}}
$$

and if $k = o(n)$, then the top order term of the numerator is $k(k-1)(k-2)n^{2k-3}$. We thus choose the parameter $c$ to be $c = \frac{k}{3}$. Then in expression (14), use the triangle inequality to estimate the deterministic part of the sum separately. It is asymptotic to $\frac{k^4}{3n^3}$, thus if $k = o(n^{3/4})$, this part of (14) goes to 0. Recall that $\mathbb{E}W = \frac{k(k-1)(k-2)}{6n^2}$, so to get a nontrivial limiting distribution, the case we are interested in is $k = \theta n^{2/3}$ for some fixed $\theta$, so limiting our considerations in a way that still allows for this case is no loss. We will now estimate each term of the non-deterministic part of (14) separately. Using the method of indicators gives

$$
\begin{aligned}
\mathbb{E}(\widetilde{M_1^2}) &= \mathbb{E}(M_1^2) - \mu_1^2 \\
&= \frac{k(n-1)n^{k-1}[(n-1)^{k-2} + (k-1)(n-2)^{k-2}] - k^2(n-1)^{2k-2}}{n^{2k-2}} \\
&= O(k).
\end{aligned}
$$

$$
\begin{aligned}
\mathbb{E}(\widetilde{M_2^2}) &= \mathbb{E}(M_2^2) - \mu_2^2 \\
&= \frac{k(k-1)(n-1)n^{k-1}[2(n-1)^{k-3} + (k-2)(k-3)(n-2)^{k-4}] - k^2(k-1)^2(n-1)^{2k-4}}{4n^{2k-2}} \\
&= O\left(\frac{k^2}{n}\right).
\end{aligned}
$$

Now to estimate (14), use Cauchy-Schwartz:

$$
\mathbb{E}\left|\frac{k(k-1)(k-2)}{6n^2} - \frac{1}{3n}[M_1 M_2 + 2M_2^2 - 2M_2]\right|
$$

$$
\begin{aligned}
&\leq O\left(\frac{k^4}{n^3}\right) + \frac{1}{3n}\mathbb{E}|\widetilde{M_1}\widetilde{M_2}| + \frac{2}{3n}\mathbb{E}\widetilde{M_2^2} + \frac{\mu_1 + 4\mu_2 - 2}{3n}\mathbb{E}|\widetilde{M_2}| + \frac{\mu_2}{3n}\mathbb{E}|\widetilde{M_1}| \\
&\leq O\left(\frac{k^4}{n^3}\right) + \frac{1}{3n}\sqrt{\mathbb{E}\widetilde{M_1^2}\mathbb{E}\widetilde{M_2^2}} + \frac{2}{3n}\mathbb{E}\widetilde{M_2^2} + \frac{\mu_1 + 4\mu_2 - 2}{3n}\sqrt{\mathbb{E}\widetilde{M_2^2}} + \frac{\mu_2}{3n}\sqrt{\mathbb{E}\widetilde{M_1^2}} \\
&\leq O\left(\frac{k^4}{n^3}\right) + O\left(\frac{k^{3/2}}{n^{3/2}}\right) + O\left(\frac{k^2}{n^2}\right) + O\left(\frac{k^2}{n^{3/2}}\right) + O\left(\frac{k^{3/2}}{n^2}\right).
\end{aligned}
$$

In particular, if $k = \theta n^{2/3}$, this expression goes to 0 as $n \to \infty$.

The same kind of analysis is used to bound (15). Again, recentering all the random variables about their means and using our choice of $c$, the expression (15) turns into:

$$
\frac{1}{n}[\widetilde{M_0}\widetilde{M_3} + \mu_3\widetilde{M_0} + \mu_0\widetilde{M_3} + \widetilde{M_1}\widetilde{M_3} + \mu_3\widetilde{M_1} + \mu_1\widetilde{M_3}] - \sum_{\ell=3}^{k}\binom{\ell}{3}\widetilde{M_\ell} + \frac{\mu_0\mu_3 + \mu_1\mu_3}{n} - \frac{k(k-1)(k-2)}{6n^2}. \tag{16}
$$

The deterministic part approaches 0:

$$
\begin{aligned}
\frac{\mu_0\mu_3 + \mu_1\mu_3}{n} - \frac{k(k-1)(k-2)}{6n^2} &= \frac{(n-1)^{2k-4}k(k-1)(k-2)(n-1+k)}{6n^{2k-1}} - \frac{k(k-1)(k-2)}{6n^2} \\
&= O\left(\frac{k^4}{n^3}\right).
\end{aligned}
$$

In order to estimate the rest, the following moments are needed:

$$
\begin{aligned}
\mathbb{E}(\widetilde{M}_0^2) &= \mathbb{E}(M_0^2) - \mu_0^2 \\
&= \frac{(n-1)^k n^{k-1} + (n-2)^k (n-1) n^{k-1} - (n-1)^{2k}}{n^{2k-2}} \\
&= O(n).
\end{aligned}
$$

$$
\begin{aligned}
\mathbb{E}(\widetilde{M}_1^2) &= \mathbb{E}(M_1^2) - \mu_1^2 \\
&= \frac{n^{k-1}[k(n-1)^{k-1} + k(k-1)(n-1)(n-2)^{k-2}] - k^2(n-1)^{2k-2}}{n^{2k-2}} \\
&= O(k).
\end{aligned}
$$

$$
\begin{aligned}
\mathbb{E}(\widetilde{M}_2^2) &= \mathbb{E}(M_2^2) - \mu_2^2 \\
&= \frac{k(k-1)(n-1)n^{k-1}[2(n-1)^{k-3} + (k-2)(k-3)(n-2)^{k-4}] - k^2(k-1)^2(n-1)^{2k-4}}{4n^{2k-2}} \\
&= O\left(\frac{k^2}{n}\right).
\end{aligned}
$$

$$
\begin{aligned}
\mathbb{E}(\widetilde{M}_3^2) &= \mathbb{E}(M_3^2) - \mu_3^2 \\
&= \frac{n^{k-1}k(k-1)(k-2)[6(n-1)^{k-3} + (k-3)(k-4)(k-5)(n-1)(n-2)^{k-6}] - k^2(k-1)^2(k-2)^2(n-1)^{2k-6}}{36n^{2k-2}} \\
&= O\left(\frac{k^3}{n^2}\right).
\end{aligned}
$$

$$
\begin{aligned}
\mathbb{E}(\widetilde{M}_\ell^2) &= \mathbb{E}(M_\ell^2) - \mu_\ell^2 \\
&= \frac{n^{k-1}[(n-1)^{k-\ell}\binom{k}{\ell} + \binom{k}{\ell}\binom{k-\ell}{\ell}(n-1)(n-2)^{k-2\ell}] - \binom{k}{\ell}^2(n-1)^{2k-2\ell}}{n^{2k-2}} \\
&= O\left(\frac{k^\ell}{n^{\ell-1}}\right) \qquad \text{(for } 4 \le \ell \le \tfrac{k}{2}\text{)}.
\end{aligned}
$$

$$
\begin{aligned}
\mathbb{E}(\widetilde{M}_\ell^2) &= \mathbb{E}(M_\ell^2) - \mu_\ell^2 \\
&= \frac{(n-1)^{k-\ell}n^{k-1}\binom{k}{\ell} - \binom{k}{\ell}^2(n-1)^{2k-2\ell}}{n^{2k-2}} \\
&= O\left(\frac{k^\ell}{n^{\ell-1}}\right) \qquad \text{(for } \tfrac{k}{2} < \ell \le k\text{)}.
\end{aligned}
$$

Looking back at (16) and using the triangle inequality and Cauchy-Schwarz on most of the terms, the following go to zero: $\mathbb{E}|\frac{1}{n}\widetilde{M}_3\widetilde{M}_0|$, $\mathbb{E}|\frac{1}{n}\mu_3\widetilde{M}_0|$, $\mathbb{E}|\frac{1}{n}\widetilde{M}_1\widetilde{M}_3|$, $\mathbb{E}|\frac{1}{n}\mu_3\widetilde{M}_1|$, and $\mathbb{E}|\frac{1}{n}\mu_1\widetilde{M}_3|$. We have already dealt with the deterministic part. Next, consider the term

$$
\mathbb{E}\left|\frac{1}{n}\mu_0\widetilde{M}_3 - \widetilde{M}_3\right| \tag{17}
$$

which is the remaining $\mu_0 \widetilde{M}_3$ term of (16) and the first summand of $\sum_{\ell=3}^{k} \binom{\ell}{3} \widetilde{M}_\ell$.

$$\mathbb{E}\left|\frac{1}{n}\mu_0\widetilde{M}_3 - \widetilde{M}_3\right| = \mathbb{E}\left|\left(\frac{(n-1)^k}{n^k} - 1\right)\widetilde{M}_3\right|$$
$$\leq O\left(\frac{k^{3/2}}{n^2}\right).$$

Finally, we have to deal with $\mathbb{E}\sum_{\ell=4}^{k}\binom{\ell}{3}|\widetilde{M}_\ell|$. Treat this expression in two parts. No matter what $n$ is, estimate the first 24 summands by $\binom{27}{3} \cdot O\left(\frac{k^2}{n^{3/2}}\right)$. The rest is bounded as follows:

$$\sum_{\ell=28}^{k}\binom{\ell}{3}\mathbb{E}|\widetilde{M}_\ell| \leq k \cdot k^3 \cdot O\left(\frac{k^{14}}{n^{27/2}}\right)$$
$$= O\left(\frac{k^{18}}{n^{27/2}}\right)$$

so if $k = o\left(n^{3/4}\right)$ as above, then this goes to 0. This proves:

**Theorem 8.** *Let $k$ balls be dropped independently and uniformly into $n$ boxes, where $W$ is the number of triples of balls sharing a box. Then,*

$$\left\|\mathcal{L}(W) - Poi_{\left(\binom{k}{3}\frac{1}{n^2}\right)}\right\|_{TV} = O\left(\frac{k^4}{n^3}\right).$$

*In particular, if $k = \theta n^{2/3}$, then $W$ has a non-trivial limiting Poisson distribution.*

**Remark:** It is instructive to compare the bound in theorem 8 to the more general bound for multiple matches in the birthday problem that appears in theorem 6B of the book of Barbour, Holst, and Janson. First, their bound gives an explicit inequality instead of the $O\left(\frac{k^4}{n^3}\right)$ of theorem 8. More importantly, in the critical case where $k = \theta n^{2/3}$, theorem 8 gives a bound of order $n^{-1/3}$ for the error. Theorem 6B gives a bound of order $n^{-2/3}$. This suggests a more careful look at the analysis above and a healthy respect for the coupling approach.

## 6. The Coupon-Collector's Problem

In its simplest version, the coupon-collector's problem is as follows: drop $k$ balls independently and uniformly into $n$ boxes. How large should $k$ be so that there is a prescribed chance (e.g. 95%) that every box contains at least one ball? This is also an old problem: Laplace (1780) gave an effective solution. Feller [26] and David-Barton [17] present substantial developments.

Write $W$ for the number of empty boxes, and use $\theta$ to write $k = n \log n + \theta n$. We introduce the new notation $k_i$ for the number of balls in box $i$, and use $N_i$ to denote the number of boxes with $i$ balls. Note that $\mathbb{E}(W) = n(1 - n^{-1})^k \sim e^{-\theta}$. Hence we take $\lambda = e^{-\theta}$. Make an exchangeable pair by choosing a ball at random, choosing a box at random, and putting the chosen ball into the chosen box. As before, compute:

$$\mathbb{P}(W' = W - 1|\omega) = \frac{(k - N_1)W}{kn}$$

and

$$\mathbb{P}(W' = W + 1|\omega) = \frac{N_1(n - W - 1)}{kn}$$

Taking $c = n$ gives the error for Poisson approximation as

$$\sup\{\mathbb{E}(Af(W + 1) - Bf(W)) : 0 \le f \le 1\}$$

where

$$A = e^{-\theta} - \frac{N_1(n - W - 1)}{k}$$

and

$$B = W - \frac{(k - N_1)W}{k} = \frac{N_1 W}{k}.$$

By the method of indicators one can compute $\mathbb{E}N_1 W = \sum_{i \ne j} \mathbb{P}(k_i = 1 \wedge k_j = 0)$, where $k_i$ is the number of balls in box $i$. Hence,

$$
\begin{aligned}
\mathbb{E}|B| &= \frac{1}{k} n(n - 1) \binom{k}{1} \frac{1}{n} \left(1 - \frac{2}{n}\right)^{k-1} \\
&\le n(1 - 2n^{-1})^{k-1} \\
&\le \frac{e^{-2\theta}}{n - 2}.
\end{aligned}
$$

Now $\mathbb{E}|A| \le \mathbb{E}|C| + \mathbb{E}|D|$, where

$$C = e^{-\theta} - \frac{(e^{-\theta} \log n)(n - W - 1)}{k}$$

and

$$D = \frac{(e^{-\theta} \log n)(n - W - 1)}{k} - \frac{N_1(n - W - 1)}{k}.$$

We have

$$C = \frac{e^{-\theta}(n \log n + \theta n - n \log n + (W + 1) \log n)}{k} = \frac{e^{-\theta}(\theta n + (W + 1) \log n)}{k}$$

and $E(W) \leq e^{-\theta}$. Hence

$$\mathbb{E}|C| \leq \frac{e^{-\theta}(\theta n + (e^{-\theta} + 1)\log n)}{n\log n + \theta n} \sim \frac{\theta e^{-\theta}}{\log n}.$$

So what remains is to bound $\mathbb{E}|D|$. Observe that $\mathbb{E}|D| \leq \frac{n\log n}{k}\mathbb{E}\left|e^{-\theta} - \frac{N_1}{\log n}\right|$.

$$\begin{aligned}
\left|\mathbb{E}\left(\frac{N_1}{\log n}\right) - e^{-\theta}\right| &= \left|\left(n + \frac{\theta n}{\log n}\right)\left(1 - \frac{1}{n}\right)^{k-1} - e^{-\theta}\right| \\
&\leq \frac{n\theta e^{-\theta}}{(n-1)\log n} + \left|n\left(1 - \frac{1}{n}\right)^{k-1} - e^{-\theta}\right| \\
&\leq \frac{n\theta e^{-\theta}}{(n-1)\log n} + \frac{e^{-\theta}}{n-1} + |e^{\alpha} - e^{-\theta}|
\end{aligned}$$

where $\alpha = \log n + k\log(1 - n^{-1})$. Now

$$\begin{aligned}
\alpha &= \log n + (n\log n + \theta n)\log\left(1 - \frac{1}{n}\right) \\
&= \log n + (n\log n + \theta n)\left[-\frac{1}{n} - \frac{1}{2n^2} - \frac{1}{3n^3} - \cdots\right] \\
&= -\theta + (n\log n + \theta n)\left[-\frac{1}{2n^2} - \frac{1}{3n^3} - \cdots\right]
\end{aligned}$$

In particular, $\alpha \leq -\theta$. Hence

$$\begin{aligned}
|e^{\alpha} - e^{-\theta}| &\leq e^{\max\{\alpha, -\theta\}}|\alpha - (-\theta)| \\
&\leq e^{-\theta}(n\log n + \theta n)\left[\frac{1}{2n^2} + \frac{1}{3n^3} + \cdots\right] \\
&\leq \frac{e^{-\theta}(\log n + \theta)}{2(n-1)}.
\end{aligned}$$

So if $\frac{N_1}{\log n}$ concentrates at its expectation, there is a bound on the error term. This is seen as follows:

$$\mathbb{E}\left|\frac{N_1}{\log n} - \mathbb{E}\left(\frac{N_1}{\log n}\right)\right| \leq \frac{\sqrt{\text{Var}(N_1)}}{\log n}.$$

Now $N_1 = \sum_{i=1}^{n}\xi_i$, where $\xi_i = \mathbb{I}\{k_i = 1\}$, $k_i$ being the number of balls in box $i$. So

$$\text{Var}(N_1) = n\text{Var}(\xi_1) + n(n-1)\text{Cov}(\xi_1, \xi_2).$$

Let

$$p = \mathbb{P}(\xi_1 = 1) = \binom{k}{1}\frac{1}{n}\left(1 - \frac{1}{n}\right)^{k-1} \sim e^{-\theta}\log n$$

and

$$\rho = \mathbb{P}(\xi_1 = 1 \wedge \xi_2 = 1) = k(k-1)\frac{1}{n^2}\left(1 - \frac{2}{n}\right)^{k-2}.$$

Then $\mathrm{Var}(\xi_1) = np(1-p) \leq np$ and $\mathrm{Cov}(\xi_1, \xi_2) = n(n-1)(\rho - p^2)$. Now

$$
\begin{aligned}
\rho - p^2 &= \frac{k(k-1)}{n^2}\left(1 - \frac{2}{n}\right)^{k-2} - \frac{k^2}{n^2}\left(1 - \frac{2}{n} + \frac{1}{n^2}\right)^{k-1} \\
&\leq \frac{k(k-1)}{n^2}\left(1 - \frac{2}{n}\right)^{k-2} - \frac{k^2}{n^2}\left(1 - \frac{2}{n}\right)^{k-1} \\
&= \left(1 - \frac{2}{n}\right)^{k-2}\left[\frac{k^2 - k}{n^2} - \frac{k^2}{n^2}\left(1 - \frac{2}{n}\right)\right] \\
&\leq \left(1 - \frac{2}{n}\right)^{k-2}\frac{2k^2}{n^3}
\end{aligned}
$$

Hence

$$
\begin{aligned}
\mathrm{Var}(N_1) &\leq k\left(1 - \frac{1}{n}\right)^{k-1} + \frac{2k^2}{n}\left(1 - \frac{2}{n}\right)^{k-2} \\
&\leq \frac{n(\log n + \theta)e^{-\theta}}{n-1} + \frac{2(\log n + \theta)^2 e^{-2\theta}}{n-2} \\
&\sim e^{-\theta}\log n.
\end{aligned}
$$

Putting everything together shows

**Theorem 9.** *Drop $k$ balls independently and uniformly into $n$ boxes. If $k = n\log n + \theta n$, and $W$ denotes the number of empty boxes, then*

$$\|\mathcal{L}(W) - Poi_{e^{-\theta}}\|_{TV} \leq O\left(\frac{e^{-\theta}}{\sqrt{\log n}}\right).$$

**Remark:** The analysis above can be sharpened to give an error of order $\frac{\log n}{n}$ in theorem 9. See the remarks to example 2 in section 9.2. It is instructive to compare this bound with the results of theorem 6D in the book by Barbour, Holst, and Janson. In the critical case where $k = n\log n + \theta n$, theorem 6D gives an explicit inequality which leads to an error of order $\frac{\log n}{n}$ as well.

## 7. Multivariate Poisson Approximation

This section gives useful bounds on the approximation of a vector of integer-valued random variables by a vector with independent Poisson coordinates. As an example, we treat the joint distribution of the number of fixed points and the number of times $i$ goes to $i+1$ in a random permutation.

**Proposition 10.** *Let $W = (W_1, \ldots, W_d)$ be a random vector with values in $\mathbb{N}^d$ and $\mathbb{E}(W_i) = \lambda_i < \infty$. Let $Z = (Z_1, \ldots, Z_d)$ have independent coordinates with $Z_i \sim Poi_{\lambda_i}$. Let $W' = (W'_1, \ldots, W'_d)$ be defined on the same probability space as $W$ with $(W, W')$ an exchangeable pair. Then for all $h : \mathbb{N}^d \to [0, 1]$,*

$$|\mathbb{E}h(W) - \mathbb{E}h(Z)| \leq \sum_{k=1}^d \alpha_k \left[ \mathbb{E}\,|\lambda_k - c_k \mathbb{P}_W(A_k)| + \mathbb{E}\,|W_k - c_k \mathbb{P}_W(B_k)| \right] \quad (18)$$

*with $\alpha_k = \min(1, 1.4\lambda_k^{-1/2})$, any choice of the $\{c_k\}$, and*

$$\begin{aligned} A_k &= \{W'_k = W_k + 1, W_j = W'_j \text{ for } k+1 \leq j \leq d\} \\ B_k &= \{W'_k = W_k - 1, W_j = W'_j \text{ for } k+1 \leq j \leq d\}. \end{aligned}$$

*Proof.* Without loss of generality, assume that the $Z$'s and the $W$'s are defined on the same space and are independent. For $h : \mathbb{N}^d \to [0, 1]$, define $f_1, f_2, \ldots, f_d$ on $\mathbb{N}^d$ as follows. Fix $k \in \{1, \ldots, d\}$ and fix $w_1, \ldots, w_{k-1}, w_{k+1}, \ldots, w_d \in \mathbb{N}$. Define $f_k(w_1, \ldots, w_{k-1}, 0, w_{k+1}, \ldots, w_d) = 0$ and

$$\begin{aligned} \lambda_k f_k(w_1, \ldots, w_{k-1}, n, w_{k+1}, \ldots, w_d) &= (n-1)f_k(w_1, \ldots, w_{k-1}, n-1, w_{k+1}, \ldots, w_d) \\ &\quad + h(w_1, \ldots, w_{k-1}, n-1, w_{k+1}, \ldots, w_d) \quad (19) \\ &\quad - \mathbb{E}h(w_1, \ldots, w_{k-1}, Z_k, w_{k+1}, \ldots, w_d). \end{aligned}$$

Lemma (1) shows that $\|f_k\|_\infty \leq \min(1, 1.4\lambda_k^{1/2})$.

Now express the difference $\mathbb{E}h(W) - \mathbb{E}h(Z)$ as a telescoping sum:

$$\mathbb{E}h(W) - \mathbb{E}h(Z) = \sum_{k=1}^d [\mathbb{E}h(Z_1, \ldots, Z_{k-1}, W_k, \ldots, W_d) - \mathbb{E}h(Z_1, \ldots, Z_k, W_{k+1}, \ldots, W_d)]. \quad (20)$$

From the definition of $f_k$ given in (19),

$$\mathbb{E}h(Z_1, \ldots, Z_{k-1}, W_k, \ldots, W_d) - \mathbb{E}h(Z_1, \ldots, Z_k, W_{k+1}, \ldots, W_d)$$

$$= \mathbb{E}[\lambda_k f_k(Z_1, \ldots, Z_{k-1}, W_k + 1, \ldots, W_d) - W_k f_k(Z_1, \ldots, Z_{k-1}, W_k, \ldots, W_d)].$$

Now, for $A_k$ and $B_k$ as defined above, define

$$T_k = c_k[f_k(Z_1, \ldots, Z_{k-1}, W'_k, \ldots, W'_d)\delta_{A_k} - f_k(Z_1, \ldots, Z_{k-1}, W_k, \ldots, W_d)\delta_{B_k}]$$

for some $c_k$, and note that by antisymmetry, $\mathbb{E}(T_k) = 0$. If $\mathcal{F}$ denotes the $\sigma$-field generated by $Z_1, \ldots, Z_d, W_1, \ldots, W_d$, then

$$\begin{aligned} \mathbb{E}(T_k|\mathcal{F}) &= c_k[f_k(Z_1, \ldots, Z_{k-1}, W_k + 1, W_{k+1}, \ldots, W_d)\mathbb{P}_W(A_k) \quad (21) \\ &\quad - f_k(Z_1, \ldots, Z_{k-1}, W_k, \ldots, W_d)\mathbb{P}_W(B_k)]. \end{aligned}$$

As $\mathbb{E}(T_k) = 0$, the expression (21) can be inserted into the summand in (20), yielding

$$\mathbb{E}h(W) - \mathbb{E}h(Z) = \sum_{k=1}^{d} \mathbb{E}[(\lambda_k - c_k\mathbb{P}(A_k))f_k(Z_1,\ldots,Z_{k-1},W_k+1,W_{k+1},\ldots,W_d) \quad (22)$$
$$-(W_k - c_k\mathbb{P}(B_k))f_k(Z_1,\ldots,Z_{k-1},W_k,\ldots,W_d)].$$

From (22), the claim (18) follows from lemma (1). □

**Example (d=2):** Choose $\sigma$ uniformly in $S_n$; let $W_1 = |\{i : \sigma(i) = i\}|$ and $W_2 = |\{i : \sigma(i) = i+1\}|$. For $W_2$, count cyclically, thus $\sigma(n) = 1$ contributes to $W_2$. It is intuitively clear that $W_1$ and $W_2$ have the same marginal distributions and that they are asymptotically independent. Proposition (10) allows a proof with error terms. Note that $\mathbb{E}(W_1) = \mathbb{E}(W_2) = 1$.

Create an exchangeable pair $(W_1', W_2')$ and $(W_1, W_2)$ by following $\sigma$ with a randomly chosen transposition. The sets $A_k$, $B_k$ are

$$
\begin{aligned}
A_1 &= \{W_1' = W_1 + 1, W_2' = W_2\} \\
B_1 &= \{W_1' = W_1 - 1, W_2' = W_2\} \\
A_2 &= \{W_2' = W_2 + 1\} \\
B_2 &= \{W_2' = W_2 - 1\}.
\end{aligned}
$$

Then

$$\mathbb{P}(A_1|\sigma) = \frac{2(n - W_1 - W_2 - W_3)}{n(n-1)},$$

where

$$W_3 = |\{i : \sigma^2(i) = i, \sigma(i) \neq i\} \cup \{i : \sigma(i) = i+1, \sigma(i+1) \neq i\} \cup \{i : \sigma(i) = \sigma^{-1}(i)+1\}|,$$

and

$$\mathbb{P}(B_1|\sigma) = \frac{2W_1(n-1-W_1-W_2)}{n(n-1)} + \frac{2W_4}{n(n-1)},$$

where

$$W_4 = |\{i : \sigma(i) = i, \sigma(i-1) = i-1\}|.$$

The calculations for $\mathbb{P}(A_2|\sigma)$ and $\mathbb{P}(B_2|\sigma)$ are analogous to those carried out in section 4. Choosing $c_k = \frac{n-1}{2}$ for $k = 1, 2$ in proposition (10) then yields the estimate

$$\|\mathcal{L}(W_1, W_2) - \mathcal{L}(Z_1, Z_2)\|_{TV} \leq \frac{13}{n},$$

where $(Z_1, Z_2)$ are independent Poisson random variables, each with mean 1.

## 8. Stein's Method for Poisson Process Approximation

### 8.1. Introduction

We begin with an example: consider, as in section (4), a random permutation $\sigma$ on $n$ letters. For $1 \leq i \leq n$, let $X_i = \begin{cases} 1 & \text{if } \sigma(i) = i, \\ 0 & \text{otherwise,} \end{cases}$ and let $W_j = \sum_{i=1}^{j} X_i$; i.e., $W$ is the number of matches up to time $j$. In the same way as in section (4), one can show that $W_j$ has an approximate Poisson distribution with parameter $\frac{j}{n}$. In this section we show that the point process $Y_t = W_{\lfloor \frac{t}{n} \rfloor}$ for $0 \leq t \leq 1$ converges to a Poisson process of rate 1. Similar results can be derived for the birthday problem: the number of birthday matches with $\lfloor t\sqrt{n} \rfloor$ people converges to a Poisson process of rate 1. For the coupon-collector's problem, the number of new coupons collected from time $n \log n$ to time $n \log n + tn$ converges to a Poisson process of rate $e^{-t}$.

Recall that if $\mathcal{Y}$ is a complete separable metric space and $\mu$ is a measure on $\mathcal{Y}$ which is finite on bounded sets, a Poisson process of rate $\mu$ is a random discrete measure on $\mathcal{Y}$ with the following properties:

- The number of points $N(A)$ in a bounded set $A$ has a Poisson distribution with parameter $\mu(A)$.
- If $A_1, \ldots, A_k$ are disjoint bounded sets, then $\{N(A_i)\}_{i=1}^{k}$ are independent.

Useful introductions to Poisson processes are found in Kingman [36] and Resnick [44]. The first use of Stein's method in Poisson process approximation was in 1988 in the paper [5] by Barbour, and the last chapter of Barbour, Holst, and Janson [8] has a wonderful development of Stein's method for Poisson process approximation, using dependency graphs and the coupling approach. We show here how the method of exchangeable pairs can be used. Again, there are other approaches to proving Poisson process approximation. We find the papers [39] and [40] by Kurtz particularly helpful. For a connection with Palm theory, see [15].

Returning to our example, once we have the Poisson process approximation in place, one has various 'off the shelf' results, e.g.

- If there are $k$ matches, their location is distributed as $k$ random points in $[0, 1]$.
- One has the limiting distribution of the minimum (or maximum) distance between two matching times (see Gilbert and Pollak [29], Feller [26]).

Furthermore, one can harness some standard machinery for transformations of Poisson processes. Perhaps the three most important such results are the following: Let a Poisson process of rate $\mu$ be given on $\mathcal{Y}$.

- **(Mapping)** If $T : \mathcal{Y} \to \mathcal{X}$ is a proper map then the image process is Poisson on $\mathcal{X}$ with rate $\mu^{T^{-1}}$.

- **(Spraying)** Let $K(y, dx)$ be a stochastic kernel on $\mathcal{Y} \times \mathcal{F}_{\mathcal{X}}$. Replace each point $y$ in the $\mathcal{Y}$ process by a point chosen from $K(y, dx)$. The new process is Poisson with rate $\int K(y, dx)\mu(dy)$.

- **(Thinning)** Delete points in the $\mathcal{Y}$ process with probability $P(y)$. The new process is Poisson with rate $\int_A P(y)\mu(dy)$.

With these examples as our motivation, we turn to a rigorous account. Let $\mathcal{I}$ be a finite index set, and let $\{X_i\}_{i \in \mathcal{I}}$ be binary indicators with some joint distribution and with $\mathbb{P}(X_i = 1) = p_i$. Let $\{W_i\}_{i \in \mathcal{I}}$ be independent Poisson random variables with parameters $p_i$. We give a development which allows us to bound the distance

$$\|\mathcal{L}(X_i : i \in \mathcal{I}) - \mathcal{L}(W_i : i \in \mathcal{I})\|_{TV}$$

between the processes. This depends on the construction of an exchangeable pair $(X, X')$. In the matching problem, X is constructed from a random permutation $\sigma$ and $X'$ is constructed by making a random transposition in $\sigma$.

### *8.2. The basic set-up*

Let $\mathcal{I}$ be a finite set, and let $\xi = \sum_{i \in \mathcal{I}} x_i \delta_i$ be a configuration on $\mathcal{I}$. Thus $\xi$ is a measure on $\mathcal{I}$ putting mass $x_i$ at $i$, with $x_i \in \mathbb{N}$. In fact, in our applications we will typically consider the case $x_i \in \{0, 1\}$. Let $\mathcal{X}_o$ be the bounded measurable functions on configurations. Thus if $f : \mathbb{N} \to \mathbb{R}$ is given, then $\xi \mapsto \sum f(i)x_i = \int f d\xi$ is in $\mathcal{X}_o$, as are $\max\{i : x_i > 0\}$ and $\int f(i, j)d\xi d\xi$ for $f : \mathbb{N} \times \mathbb{N} \to \mathbb{R}$.

Let $\mathbb{E}_o$ be the Poisson expectation operator on $\mathcal{X}_o$. In other words, consider a probability space $(\Omega, \mathcal{A}, \mathbb{P})$ and an induced probability on the space of configurations on $\mathcal{I}$ given by writing $\xi = \sum x_i(\omega)\delta_i$ where the $x_i \sim Poi(p_i)$. Let $\mathcal{X}$ be the space of bounded measurable functions on $\Omega$ with corresponding expectation operator $\mathbb{E}$. As in section 2, there will be a map $\beta : \mathcal{X}_o \to \mathcal{X}$ assigning a random variable to a function on configurations:

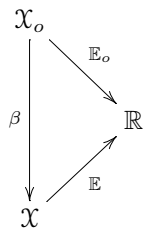

In the matching example, the probability space is $S_n$ with the uniform distribution, and the association $\sigma \longleftrightarrow \xi_\sigma = \sum x_i(\sigma)\delta_i$ gives rise to the map $\beta T(\sigma) = T(\xi_\sigma)$ for $T \in \mathcal{X}_o$. The aim is to show that $\mathbb{E}\beta \doteq \mathbb{E}_o$. To do this we

introduce

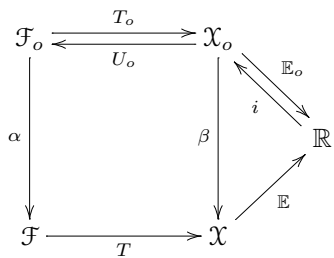

with $T_o$ the Stein operator

$$T_o h(\xi) = \sum_{i \in \mathcal{I}} p_i [h(\xi + \delta_i) - h(\xi)] + \sum_{i \in \mathcal{I}} x_i [h(\xi - \delta_i) - h(\xi)] \qquad (23)$$

and $\mathcal{F}_o$ the subset of elements $h \in \mathcal{X}_o$ with $T_o h$ bounded. As explained in Barbour, Holst, and Janson, this $T_o$ is the generator of an immigration and death process $Z(t)$ on $\mathcal{I}$ where particles are born at site $i$ at exponential rate $p_i$ and die at rate one. This process has the Poisson process of rate $\{p_i\}$ as a stationary distribution. This fact (or a straight-forward direct computation) shows that $\mathbb{E}_o T_o = 0$.

Note that for $|\mathcal{I}| = 1$ the operator is slightly different from the usual Stein operator for a Poisson process. This $T_o$ is a "second difference" operator. We use it because of already available bounds on its inverse, but also as an introduction to this line of development. The inverse can be given a stochastic representation

$$U_o f(\xi) = - \int_0^\infty [\mathbb{E}(f(Z(t)) | \xi) - \mathbb{E}_o f] dt$$

where $Z(t)$ is the immigration and death process started at $\xi$. Barbour, Holst, and Janson (pg. 209) show that

$$T_o U_o f(\xi) = f(\xi) - \mathbb{E}_o f.$$

They also show that if $f(\xi) = \delta_A(\xi)$ for some set of configurations $A$, then

$$\sup\{|U_o f(\xi + \delta_i) - U_o f(\xi)| : i \in \mathcal{I}, \xi \text{ a configuration}\} \le 1.$$

As before, let $\mathcal{F}$ be the bounded anti-symmetric functions on $\Omega \times \Omega$, and construct $T$ via a symmetric probability $\mathbb{Q}$ with margins $\mathbb{P}$:

$$T f(\omega) = \mathbb{E}_{\mathbb{Q}}(f(\omega, \omega') | \omega).$$

We may choose $\alpha$ to be any linear map $\alpha : \mathcal{F}_o \to \mathcal{F}$. Stein's identity is

$$\mathbb{E}(h(\xi)) = \mathbb{E}_o(h(\xi)) + \mathbb{E}[T\alpha - \beta T_o] U_o h.$$

To proceed, we need to make an intelligent choice of $\alpha$, make explicit the difference $T\alpha - \beta T_o$, and use this expression to bound the remainder term.

Because of the different form of the Stein operator, the choice of $\alpha$ which seems most useful is:

$$\alpha f(\omega, \omega') = c[\beta f(\xi_\omega) - \beta f(\xi_{\omega'})]\delta(\omega, \omega') - c[\beta f(\xi_{\omega'}) - \beta f(\xi_\omega)]\delta(\omega', \omega)$$

where

$$\delta(\omega, \omega') = \begin{cases} 1 & \text{if } \xi_{\omega'} = \xi_\omega + \delta_i \text{ for some } i, \\ 0 & \text{otherwise.} \end{cases}$$

Thus,

$$\begin{aligned} \beta T_o h(\omega) &= \sum_{i\in\mathcal{I}} p_i[h(\xi_\omega + \delta_i) - h(\xi_\omega)] - \sum_{i\in\mathcal{I}} x_i[h(\xi_\omega - \delta_i) - h(\xi_\omega)] \\ T\alpha h(\omega) &= c\sum_{i\in\mathcal{I}}[h(\xi_\omega + \delta_i) - h(\xi_\omega)]\mathbb{Q}(\xi_{\omega'} = \xi_\omega + \delta_i \,|\, \omega) \\ &\quad -c\sum_{i\in\mathcal{I}}[h(\xi_\omega) - h(\xi_\omega - \delta_i)]\mathbb{Q}(\xi_{\omega'} = \xi_\omega - \delta_i \,|\, \omega). \end{aligned}$$

In our application, $h = U_o f$, thus $h$ has bounded differences. It then follows that if $c$ and the probability $\mathbb{Q}$ can be chosen so that

$$\begin{aligned} |c\mathbb{Q}(\xi_{\omega'} = \xi_\omega + \delta_i \,|\, \omega) - p_i| &< \epsilon_i \quad \text{and} \\ |c\mathbb{Q}(\xi_{\omega'} = \xi_\omega - \delta_i \,|\, \omega) - x_i| &< \epsilon'_i \end{aligned}$$

then the error in the Poisson approximation is bounded by $\sum_{i\in\mathcal{I}}(\epsilon_i + \epsilon'_i)$.

Of course, in any application, some work is required. Here is a full version of the matching problem.

**Theorem 11.** *Let $X_i(\sigma), 1 \leq i \leq n$ indicate whether or not $i$ is a fixed point of the random permutation $\sigma$. Let $Y_i$ be independent Poisson random variables, each with mean $\left(\frac{1}{n}\right)$. Then*

$$\|\mathcal{L}(X_i : 1 \leq i \leq n) - \mathcal{L}(Y_i : 1 \leq i \leq n)\|_{TV} \leq \frac{4}{n}.$$

*Proof.* It is easy to see that

$$\mathbb{Q}(\xi_{\omega'} = \xi_\omega + \delta_i \,|\, \sigma) = \frac{2(1 - Z_i(\sigma))}{n(n-1)}$$

and that

$$\mathbb{Q}(\xi_{\omega'} = \xi_\omega - \delta_i \,|\, \sigma) = \frac{2X_i(\sigma)(n - W(\sigma))}{n(n-1)}$$

with $Z_i(\sigma)$ indicates whether $\sigma(\sigma(i))$ is $i$ or not, and $W = \sum X_i$. Choose $c = \frac{n-1}{2}$, and note that $p_i = 1$ for each $i$. We thus have that

$$\left| c\mathbb{Q}(\xi_{\omega'} = \xi_\omega + \delta_i \,|\, \sigma) - \frac{1}{n} \right| = \frac{Z_i(\sigma)}{n}$$

and that

$$|c\mathbb{Q}(\xi_{\omega'} = \xi_\omega - \delta_i \,|\, \sigma) - X_i(\sigma)| = \frac{X_i(\sigma)W(\sigma)}{n}.$$

It follows that if $A$ is any set of configurations,

$$
\begin{aligned}
|\mathbb{P}(\xi(\sigma) \in A) - \mathbb{E}_o(\delta_A)| &\leq \frac{1}{n}\mathbb{E}(W + 2a_2(\sigma)) + \frac{1}{n}\mathbb{E}(W^2) \\
&= \frac{4}{n}.
\end{aligned}
$$

$\square$

**Remark:** We have not treated the closely related area of compound Poisson approximation; see Erhardsson [25] and the references therein.

## 9. The coupling approach to Poisson approximation

Stein's method for Poisson approximation has been very actively developed along two lines, both quite different from the present account. The connections between the various approaches are not clear. In this section we give a brief development of the coupling method. In section 10 we give a brief development of the dependency graph method. The book-length account of Barbour, Holst, and Janson treats both coupling and dependency graphs in detail. The main reasons for including brief descriptions here are twofold: (a) as illustrated above, the many analytic tools developed by previous workers can easily be adapted for use with the present method of exchangeable pairs; (b) while the three approaches are similar, there are also real differences. Sometimes one or another approach is easier to apply or yields better bounds.

The basics of the coupling method are laid out in section 9.1. Section 9.2 treats our familiar examples: Poisson-binomial trials, the birthday problem, the coupon-collector's problem, and the matching problem. Section 9.3 illustrates one of the triumphs of the coupling approach: effortless bounds for negatively dependent summands. Our development leans heavily on lectures of Charles Stein.

### 9.1. The basics

From section 2, we know that a random variable $W$ has a Poisson distribution with parameter $\lambda$ if and only if, for each bounded measureable function $f$,

$$\mathbb{E}(Wf(W)) = \lambda\mathbb{E}(f(W + 1)).$$

This suggests that one approach to the question of whether $W$ is approximately Poisson-distributed is to compare $\mathbb{E}(Wf(W))$ and $\lambda\mathbb{E}(f(W + 1))$. The following simple lemma motivates a slightly mysterious construction to come. It appears in Stein [51].

**Lemma 12.** *Let $X = \{X_i\}_{i \in \mathcal{I}}$ be $\{0, 1\}$-valued random variables indexed by a finite set $\mathcal{I}$. Let $\lambda = \sum_{i \in \mathcal{I}} \mathbb{E}(X_i)$ and let $I$ be uniformly distributed on $\mathcal{I}$ and independent of $X$. Then, for any bounded $f$,*

$$\mathbb{E}(Wf(W)) = \lambda \mathbb{E}(f(W)|X_I = 1).$$

*Proof.*

$$
\begin{aligned}
\mathbb{E}(Wf(W)) &= |\mathcal{I}| \cdot \mathbb{E}\{\mathbb{E}(X_I|X)f(W)\} \\
&= |\mathcal{I}| \cdot \mathbb{E}\{X_I f(W)\} \\
&= |\mathcal{I}| \cdot \mathbb{E}\{X_I \mathbb{E}(f(W)|X_I = 1)\} \\
&= \lambda \mathbb{E}(f(W)|X_I = 1).
\end{aligned}
$$

$\square$

The lemma suggests comparing $W + 1$ and a random variable $W^*$ having the distribution of $W$, given that the randomly chosen coordinate $X_I$ is one. Many examples in which $W^*$ is constructed explicitly appear below. The following result gives a concrete bound on the total variation distance between the law of $W$ and $Poi_\lambda$ in terms of the difference between $W + 1$ and $W^*$.

**Proposition 13.** *Let $X$, $I$, and $W$ be defined as in lemma (12). Let $W^*$ be defined on the same probability space with*

$$\mathbb{P}(W^* = w) = \mathbb{P}(W = w|X_I = 1).$$

*Then*

$$\|\mathcal{L}(W) - Poi_\lambda\|_{TV} \leq (1 - e^{-\lambda})\mathbb{E}|W + 1 - W^*|.$$

*Proof.* Given $A \subseteq \mathbb{N} \cup \{0\}$, let $g(j) = \delta_A(j) - Poi_\lambda(A)$ and set $f = U_o g$, for $U_o$ as defined in (5). Then, from lemma (12) and equation (6),

$$
\begin{aligned}
0 &= \mathbb{E}(Wf(W) - \lambda f(W^*)) \\
&= \mathbb{E}(Wf(W) - \lambda f(W + 1) + \lambda f(W + 1) - \lambda f(W^*)) \\
&= \mathbb{P}(W \in A) - Poi_\lambda(A) + \lambda \mathbb{E}(f(W + 1) - f(W^*)).
\end{aligned}
$$

Now, for any integers $x < y$,

$$|f(x) - f(y)| \leq |f(x) - f(x+1)| + |f(x+1) - f(x+2)| + \ldots + |f(y-1) - f(y)|,$$

and from lemma (1), $|f(i+1) - f(i)| \leq \dfrac{(1 - e^{-\lambda})}{\lambda}$, thus

$$|\mathbb{P}(W \in A) - Poi_\lambda(A)| \leq (1 - e^{-\lambda})\mathbb{E}|W + 1 - W^*|.$$

$\square$

From these pieces, we get a bound for Poisson approximation by constructing a coupling $(W, W^*)$ with $W^*$ close to $W + 1$. Stein [51] and Barbour, Holst, and Janson [8] give many examples leading to good bounds in a host of problems.

We next give a brief treatment of the basic probability problems treated in previous sections.

### *9.2. Examples*

**Example 0 (Poisson-binomial trials):** Let $X_i$, $1 \leq i \leq n$ be independent with $\mathbb{P}(X_i = 1) = p_i = 1 - \mathbb{P}(X_i = 0)$. Let $\lambda = \sum p_i$ and let $W = \sum_{i=1}^{n} X_i$. Then

$$\|\mathcal{L}(W) - Poi_\lambda\|_{TV} \leq \frac{(1 - e^{-\lambda})}{\lambda} \sum_{i=1}^{n} p_i^2.$$

*Proof.* We construct $W^*$ 'backwards'. Pick $X$ from the product measure, and independently pick $I^*$ with $\mathbb{P}(I^* = i) = \frac{p_i}{\lambda}$. Set $X_{I^*}^* = 1$ and $X_i^* = X_i$ for $i \neq I^*$. We claim that

$$\mathbb{P}(I^* = i, X^* = x) = \mathbb{P}(I = i, X = x | X_I = 1).$$

Indeed, if $x_i = 0$, then both sides are 0. If $x_i = 1$, then

$$\mathbb{P}(I^* = i, X^* = x) = \frac{p_i}{\lambda} \prod_{j \neq i} p_j^{x_j} (1 - p_j)^{1 - x_j}$$

and

$$\mathbb{P}(I = i, X = x | X_I = 1) = \frac{\mathbb{P}(I = i, X = x)}{\mathbb{P}(X_I = 1)} = \frac{\frac{1}{n} p_i \prod_{j \neq i} p_j^{x_j} (1 - p_j)^{1 - x_j}}{\frac{\lambda}{n}}.$$

Note that the claim shows $\mathbb{P}(W^* = j) = \mathbb{P}(W = j | X_I = 1)$. To complete the proof, observe that $W + 1 - W^* = X_{I^*}$, and that $\mathbb{P}(X_{I^*} = 1) = \sum_i \frac{p_i^2}{\lambda}$. □

**Remarks:**

1. Observe that the bound of this example is exactly the same as the bound given by the method of exchangeable pairs.
2. The backwards construction of this example seems magical when you first encounter it. The method of exchangeable pairs seems much more straightforward (at least to us). One benefit of the coupling approach is the clean bound of proposition (13).

**Example 1 (The matching problem):** Here, let $X_i = \delta_{i\pi(i)}$ with $\pi$ a randomly chosen permutation in $S_n$. Then $W = \sum_{i=1}^{n} X_i$ is the number of fixed points in $\pi$. Construct $W^*$ by picking $I$ uniformly in $\{1, \dots, n\}$ and then setting

$$\pi^* = \begin{cases} \pi & \text{if } \pi(I) = I, \\ (IJ)\pi & \text{if } \pi(I) = J \neq I. \end{cases}$$

Then if $a_2$ is the number of transpositions in $\pi$,

$$W + 1 - W^* = \begin{cases} 1 & \text{with probability } \frac{W}{n}, \\ -1 & \text{with probability } \frac{2a_2}{n}, \\ 0 & \text{otherwise.} \end{cases}$$

Then

$$\begin{aligned} \mathbb{E}|W + 1 - W^*| &= \mathbb{E}\mathbb{E}_\pi|W + 1 - W^*| \\ &= \mathbb{E}\left[\frac{W}{n} + \frac{2a_2}{n}\right] \\ &= \frac{2}{n}, \end{aligned}$$

thus

$$\|\mathcal{L}(W) - Poi_1\|_{TV} \leq \frac{(1 - e^{-1})2}{n}.$$

**Example 2 (The coupon-collector's problem):** Here $k$ balls are dropped into $n$ boxes and $X_i$ is one or zero as box $i$ is empty or not. Then $W = \sum_i X_i$ is zero if all of the boxes are covered. To construct a coupling, pick a box uniformly and simply distribute each ball in it at random into one of the other boxes. We claim that

$$\mathbb{E}(|W + 1 - W^*|) \leq \left(1 - \frac{1}{n}\right)^k \left(1 + \frac{k}{n}\right).$$

To prove this, let $A$ be the number of balls in the box selected. Observe that $W + 1 - W^* \geq 0$ and that

- If $A = 0$, then $W + 1 - W^* = 1$.
- If $A > 0$, then $W + 1 - W^*$ is the number of the $A$ balls that fall into empty boxes.

From this,

$$\mathbb{E}(|W + 1 - W^*|) \leq \frac{\mathbb{E}(W)}{n} + \frac{\mathbb{E}(AW)}{n}.$$

Now, $\mathbb{E}(W) = n\left(1 - \frac{1}{n}\right)^k$. To compute $\mathbb{E}(AW)$, observe that by symmetry this is $\mathbb{E}(A_1 W)$, with $A_1$ the number of balls in box 1. Let $X_i$ indicate whether box $i$ is empty, thus $\sum_{i=1}^{n} X_i = W$, and let $Y_j$ indicate whether ball $j$ falls into the first box, thus $\sum_{j=1}^{k} Y_j = A_1$. Then

$$
\begin{aligned}
\mathbb{E}(A_1 W) &= \mathbb{E}(A_1 X_1) + (n-1)\mathbb{E}(A_1 X_2) \\
&= (n-1)k\mathbb{E}(Y_1 X_2) \\
&= k\left(1 - \frac{1}{n}\right)^k.
\end{aligned}
$$

Combining these bounds proves the claim.

**Remark:** If $k = n\log n + cn$ for $c$ fixed and $n$ large, then proposition (13) shows that the number of empty boxes is well approximated by a Poisson distribution with error of order $\frac{\log n}{n}$. This is better than the bound we obtained in section 6 using exchangeable pairs. In preliminary work we have shown that if the proof of section 6 is carried through by conditioning only on $W$, we get an error term of order $\frac{\log n}{n}$ as well.

**Example 3 (The birthday problem):** Here $k$ balls are dropped uniformly into $n$ boxes. Let $X_{ij}$ indicate whether balls $i$ and $j$ fall into the same box. Let $W = \sum_{i<j} X_{ij}$. To form $W^*$, pick a random pair of indices $I < J$ and move ball $J$ into the same box as ball $I$.

- If balls $I$ and $J$ are in the same box already (call this event $S$), then $|W + 1 - W^*| = 1$.
- If ball $I$ is in a box containing $a$ balls and ball $J$ is in a (different) box containing $b$ balls, then $|W + 1 - W^*| = |a - b|$.

Thus

$$
\mathbb{E}(|W + 1 - W^*|) = \mathbb{P}(S) + \mathbb{E}\left(|W + 1 - W^*|\Big| S^c\right) \cdot \mathbb{P}(S^c).
$$

Now $\mathbb{P}(S) = \sum_w \mathbb{P}(S|W=w)\mathbb{P}(W=w) = \sum_w \dfrac{w}{\binom{k}{2}}\mathbb{P}(W=w) = \dfrac{1}{\binom{k}{2}}\mathbb{E}(W) = \dfrac{1}{n}$.

For $\mathbb{E}\left(|W + 1 - W^*|\Big| S^c\right)$ we may, by symmetry, consider the absolute difference between the number of balls in the first two boxes, given that each contains at least one ball. This is just the absolute difference between the number of balls in the first and second boxes when $k - 2$ balls are dropped into $n$ boxes. By bounding this by twice the expected number of balls in the first box, namely $\frac{2k}{n}$, and bounding $\mathbb{P}(S^c)$ by 1, we get

$$
\mathbb{E}(|W + 1 - W^*|) \leq \frac{1 + 2k}{n}.
$$

**Remark:** If $k = \theta\sqrt{n}$ with $\theta$ fixed and $n$ large, then proposition (13) shows that the number of matching pairs has a Poisson distribution of mean $\frac{\theta^2}{2}$ with error of order $\frac{1}{\sqrt{n}}$.

### 9.3. Automatic coupling

One of the many contributions of the book-length treatment on Poisson approximation by Barbour, Holst, and Janson (1992) is a soft theory which leads to very neat result. Recall that a collection $\{X_i\}_{i \in \mathfrak{I}}$ of random variables is called *negatively associated* if for all disjoint $A, B \subseteq \mathfrak{I}$ and monotone functions $f, g$,

$$\mathbb{E}(f(X_i; i \in A)g(X_j; j \in B)) \leq \mathbb{E}(f(X_i; i \in A))\mathbb{E}(g(X_j; j \in B)).$$

Here, $f : \mathbb{R}^m \to \mathbb{R}$ is monotone if $x \preceq y$ implies $f(x) \leq f(y)$ where $x \preceq y$ is the coordinate order.

Two leading tools are: (a)If $\{Y_i\}$ are independent with log-concave densities, and $Z_i$ are distributed as $Y_i$ given $\sum Y_i = k$, then $\{Z_i\}$ are negatively associated. (b)Increasing functions of disjoint subsets of negatively associated variables are negatively associated.

By way of example, if $Y_i$ are independent, with $Y_i \sim Poi_{p_i\lambda}$, where $\sum_i p_i = 1$ and $N_i$ have the law of $Y_i$ given $\sum_i Y_i = k$, then $N_i$ have the law of $k$ balls dropped into boxes via a multinomial allocation. They are negatively associated. If $X_i$ is one or zero as $N_i$ is bounded above by $m$ or not (for fixed $m \geq 0$), then $X_i$ are negatively associated.

In Bose-Einstein allocation, $k$ balls are dropped into $n$ boxes so that all configurations are equally likely. If $Y_i$ are independent geometric random variables with $\mathbb{P}(Y_i = j) = qp^j$ for $0 \leq j < \infty$, then the distribution of $Y_i$ given $\sum_i Y_i = k$ has a Bose-Einstein distribution. These variables are thus negatively associated. Many further examples and tools are given in [8].

The motivation for this set-up is the following elegant theorem of Barbour, Holst, and Janson.

**Theorem 14.** *Let $\{X_i\}_{i \in \mathfrak{I}}$ be binary, negatively associated random variables. Let $W = \sum_{i \in \mathfrak{I}} X_i$ have mean $\lambda$ and variance $\sigma^2$. Then*

$$\|\mathcal{L}(W) - Poi_\lambda\|_{TV} \leq (1 - e^{-\lambda})\left(1 - \frac{\sigma^2}{\lambda}\right).$$

**Remarks:**

1. Of course, the mean equals the variance for a Poisson random variable. The theorem shows that if the mean is close to the variance, then the law is close to Poisson.

2. Mixtures of negatively associated variables are negatively associated. Diaconis and Holmes [22] apply these ideas to study the birthday and coupon-collector's problems when balls are dropped into boxes with probabilities $\{p_i\}$ and a prior is put on $p_i$.

3. For variables presented as conditional given a sum, a host of analytic techniques are available for higher order expansions and large deviations [20].

## 10. The dependency graph approach

In this section we give a brief description of the dependency graph approach to Stein's method. This method is often useful when a natural local dependence structure exists; in these cases, this approach may be easier to apply than constructing an exchangeable pair or a coupling. Following a description of the approach, an example of using a dependency graph to solve a variation of the birthday problem is given.

Let $\{X_i\}_{i=1}^n$ be a set of binary random variables, with $\mathbb{P}(X_i = 1) = p_i$. A dependency graph for $\{X_i\}$ is a graph with vertex set $I = \{1, \ldots, n\}$ and edge set $E$, such that if $I_1, I_2$ are disjoint subsets of $I$ with no edges connecting them, then $\{X_i\}_{i \in I_1}$ and $\{X_i\}_{i \in I_2}$ are independent. Let $N_i$ denote the neighborhood of $i$ in the graph; that is,

$$N_i = \{j \in I : (i, j) \in E\} \cup \{i\}.$$

This framework yields the following bound, first proved in [4]:

**Theorem 15.** *Let $\{X_i\}_{i \in I}$ be a finite collection of binary random variables with dependency graph $(I, E)$; suppose $\mathbb{P}(X_i = 1) = p_i$ and $\mathbb{P}(X_i = 1, X_j = 1) = p_{ij}$. Let $\lambda = \sum p_i$ and $W = \sum X_i$. Then*

$$\|\mathcal{L}(W) - Poi_\lambda\|_{TV} \leq \min(1, \lambda^{-1}) \left[ \sum_{i \in I} \sum_{j \in N_i \setminus \{i\}} p_{ij} + \sum_{i \in I} \sum_{j \in N_i} p_i p_j \right].$$

*Proof.* Let $A \subseteq \mathbb{N}$, and let $f = U_o \delta_A$ for $U_o$ as in section 2. Let $Z \sim Poi_\lambda$. Then

$$\begin{align} \mathbb{P}(Z \in A) - \mathbb{P}(W \in A) &= \mathbb{E}[Wf(W) - \lambda f(W+1)] \qquad (24) \\ &= \sum_{i \in I} \mathbb{E}[X_i f(W) - p_i f(W+1)]. \end{align}$$

Let $W_i = W - X_i$, and $V_i = \sum_{j \notin N_i} X_j$. Note that $X_i$ and $V_i$ are independent and that

$$X_i f(W) = X_i f(W_i + 1).$$

So from (24),

$$
\begin{aligned}
\mathbb{P}(Z \in A) - \mathbb{P}(W \in A) &= \sum_{i \in I} \mathbb{E}[(X_i - p_i) f(W_i + 1)] + p_i \mathbb{E}[f(W_i + 1) - f(W + 1)] \\
&= \sum_{i \in I} \left[ \mathbb{E}\left[ (X_i - p_i)\big(f(W_i + 1) - f(V_i + 1)\big) \right] + p_i \mathbb{E}[f(W_i + 1) - f(W + 1)] \right].
\end{aligned}
$$

By lemma (1), $|f(W_i + 1) - f(W + 1)| \leq \min(1, \lambda^{-1})$, thus

$$
|p_i \mathbb{E}[f(W_i + 1) - f(W + 1)]| \leq \min(1, \lambda^{-1}) p_i^2.
$$

(The second factor of $p_i$ occurs as $\mathbb{P}(W_i \neq W) = p_i$.) Further, writing $f(W_i + 1) - f(V_i + 1)$ as a sum of telescoping terms yields:

$$
f(W_i + 1) - f(V_i + 1) \leq \min(1, \lambda^{-1})(W_i - V_i) = \min(1, \lambda^{-1}) \left( \sum_{j \in N_i \setminus i} X_i \right),
$$

which gives

$$
|\mathbb{E}[(X_i - p_i)(f(W_i + 1) - f(V_i + 1))]| \leq \min(1, \lambda^{-1}) \sum_{j \in N_i \setminus i} (p_{ij} + p_i p_j).
$$

This proves the theorem. $\qquad\square$

**Remark:** Here we have insisted that $X_i$ be independent with $\{X_j\}_{j \notin N_i}$. The following theorem, stated and proved in Barbour, Holst, and Janson (1992), is a more general result in which it is only required that $X_i$ must not depend too strongly on $\{X_j\}_{j \notin N_i}$. The proof is very similar to that of theorem (15).

**Theorem 16.** *Let $\{X_i\}_{i \in I}$ be a finite set of binary random variables with $\mathbb{P}(X_i = 1) = p_i$, and let $\lambda = \sum p_i$ and $W = \sum X_i$. For each $i \in I$, let $I_i^w$ and $I_i^s$ be disjoint subsets of $I$ so that $(I_i^w, I_i^s, \{i\})$ partitions the set $I$. Define the random variables*

$$
\begin{aligned}
Y_i &= \sum_{j \in I_i^w} X_j \\
Z_i &= \sum_{j \in I_i^s} X_j.
\end{aligned}
$$

*Then*

$$
\|\mathcal{L}(W) - Poi_\lambda\|_{TV} \leq \min(1, \lambda^{-1}) \sum_{i \in I} [p_i^2 + p_i \mathbb{E} Z_i + \mathbb{E}(X_i Z_i)] + \min(1, \lambda^{-1}) \sum_{i \in I} \eta_i,
\tag{25}
$$

*where $\eta_i$ is any quantity satisfying*

$$
|\mathbb{E}(X_i g(Y_i + 1)) - p_i \mathbb{E} g(Y_i + 1)| \leq \eta_i
$$

*for every $g : \mathbb{N} \cup \{0\} \to [0, 1]$.*

**Remarks:**

1. The idea is to choose the sets $I_i^w$ and $I_i^s$ to be the indexes of those $X_j$ on which $X_i$ depends 'weakly' and 'strongly', respectively. If the $X_i$ are equipped with a dependency graph $(I, E)$ as is theorem (15), then using $I_i^w = I \setminus N_i$ and $I_i^s = N_i \setminus i$ recovers theorem (15) from theorem (16).

2. If $I_i^w$ is well-chosen as those $X_j$ on which $X_i$ only depends weakly, then one potentially good choice of $\eta_i$ is

$$\eta_i = \mathbb{E}|\mathbb{E}(X_i|\{X_j, j \in I_i^w\}) - p_i|.$$

**Example (A generalized birthday problem):** Let $[n] = \{1, \dots, n\}$, fix $k \geq 2$, and let $I = \{\alpha \subseteq [n] : |\alpha| = k\}$. Color the points of $[n]$ independently, choosing the color of each uniformly from a set of $c$ colors. Let

$$X_\alpha = \begin{cases} 1 & \text{if all points of } \alpha \text{ have the same color,} \\ 0 & \text{otherwise,} \end{cases}$$

and let $W = \sum_{\alpha \in I} X_\alpha$. Thus $W > 0$ means that at least one group of $k$ entries in $[n]$ all have the same color. We show that $W$ is approximately Poisson with parameter $\lambda = \binom{n}{k}c^{1-k}$.

To get the problem into the present framework, we need a dependency graph for the $X_\alpha$. The simplest choice is to observe that if $\alpha \cap \beta = \emptyset$, then $X_\alpha$ and $X_\beta$ are independent, thus one dependency graph is the one in which $N_\alpha = \{\beta : \alpha \cap \beta \neq \emptyset\}$. To apply theorem (15), note that for $\alpha \subseteq [n]$,

$$p_\alpha = \mathbb{P}(X_\alpha = 1) = c^{1-k}, \tag{26}$$

and if $|\alpha \cap \beta| = \ell \neq 0$, then

$$p_{\alpha\beta} = \mathbb{P}(X_\alpha = 1, X_\beta = 1) = c^{1-(2k-\ell)}. \tag{27}$$

Given $\alpha$, a set $\beta$ such that $|\alpha \cap \beta| = \ell$ can be chosen in $\binom{k}{\ell}\binom{n-k}{k-\ell}$ ways. Putting these observations into theorem (15) yields the following:

**Proposition 17.** *For positive integers $c$, $k$, and $n$, let $W$ be the number of monochromatic $k$-tuples if an $n$-element set is colored with $c$ colors, independently and uniformly. Let $\lambda = \binom{n}{k}c^{1-k}$. Then*

$$\|\mathcal{L}(W) - Poi_\lambda\|_{TV} \leq \min(1, \lambda^{-1}) \left[ \binom{n}{k}\sum_{\ell=1}^{k-1}\binom{k}{\ell}\binom{n-k}{k-\ell}c^{1-(2k-\ell)} + \binom{n}{k}c^{2-2k}\sum_{\ell=1}^{k}\binom{k}{\ell}\binom{n-k}{k-\ell} \right].$$

**Remarks:**

1. Restricting to the cases $k = 2, 3$, the bound of proposition (17) is of the same order as the bounds obtained in section 5.

2. The dependency graph approach seems quite well suited to this problem; theorem (17) was obtained very easily and is rather more general than the results in section 5, while those arguments were rather technical. On the other hand, the generalized matching problem presented in section 4 was not so difficult via the method of exchangeable pairs, and the dependence structure is much less clear in that problem. Although there are problems for which both approaches work well (Poisson-binomial trials are an easy example; see Barbour, Holst, and Janson for more involved examples), it seems that these two methods mostly are useful in different contexts.

## 11. Some open problems

Our treatment above has been expository and the classical problems studied were meant to help illuminate the basics of Stein's method of exchangeable pairs. A point is, this method may be used for more complex problems than the ones treated here. The following examples are all open as far as we know.

### 11.1. Variants of the matching problem

1. *Permutations with restricted positions.* Let $A$ be a fixed $n \times n$ matrix with zero-one entries. A permutation $\sigma \in S_n$ is compatible with $A$ if $\prod_{i=1}^{n} a_{i\,\sigma(i)} = 1$. If, for example, $A$ is zero on the diagonal and one elsewhere, then those $\sigma$ which are compatible with $A$ are the ones with no fixed points. The matching problem asks for the number of permutations with $k$ violations. One also may investigate this problem for more general $A$. This is a rich subject; for pointers, one may look at the classical treatment in Riordan [46] or the section on Rook polynomials in Stanley [48]. Diaconis-Graham-Holmes [21] give a survey and applications in statistics. In addition to studying classical patterns such as the menage problem, one may ask for necessary and/or sufficient conditions on $A$ to guarantee a Poisson limit.

2. *Other actions.* Let $G$ be a finite group and $X$ a finite set with $G$ acting on $X$. One may mimic the development in the example above by considering a matrix $A$ indexed by $X$ (i.e., with entries $a_{x\,y}$) with zero-one entries. Say $g \in G$ agrees with $A$ if $\prod_{x} a_{x\,gx} = 1$. One could choose $g$ uniformly in $G$ and ask for the distribution of violations. An example of statistical interest (Daniels's test – see [34]) has $G$ the symmetric group $S_n$ and $X$ is the set of ordered pairs $\{(i,j) : i \neq j\}$ with $g(i,j) = (g(i), g(j))$.

3. *Non-uniform permutations.* All of the development above has been for uniformly chosen permutations (or group elements). A nice problem is to develop Poisson approximation for natural families of non-uniform distributions. Here is a specific example: Diaconis [19] and Critchlow [16] discuss a variety of metrics on the permutation group. For example, one may take Spearman's $\rho$, which boils down to the squared differences $\sum (\sigma(i) - \sigma'(i))^2$. Fixing a metric $d$, one can define a probability distribution on permutations by

$$P_\theta(\sigma) = z^{-1} \theta^{d(\sigma, \sigma_o)}, \qquad 0 < \theta < 1,$$

where $z$ is a normalizing constant and $\sigma_o$ is the center of the distribution. These Mallows models have been widely applied; for extensive references, see Fligner and Verducci [27] or Diaconis and Ram [24].

The problem is this: for $\theta$, $\sigma_o$, and $d$ fixed, pick $\sigma$ from $\mathbb{P}_\theta$ and look at the distribution of the number of fixed points. We conjecture that this will be Poisson under very mild conditions. To use Stein's method, an exchangeable pair must be formed. One simple way to do this is to pick $\sigma$ from $\mathbb{P}_\theta$ and then make a random transposition $\tau$. Use the Metropolis algorithm to get the second coordinate distribution to be $\mathbb{P}_\theta$ as well. Alternatively, the Gibbs sampler may be used. Of course, the problems outlined in the first two items above may be studied under these non-uniform models. Diaconis and Holmes [22] give some examples and further motivation.

### 11.2. Variants of the coupon-collector's problem

1. *Non-uniform allocations.* One can drop $k$ balls into $n$ boxes with the probability of dropping a ball into box $i$ equal to $p_i$. This is studied by classical arguments (and without error terms) in Rosén [47]. Diaconis and Holmes [22] study Bose-Einstein allocation. They also mention the problem of putting a prior on $(p_1, \ldots, p_n)$ and studying the Bayes distribution.

2. *Complexes.* The balls may be dropped in $j$ at a time such that within each group of balls, all must land in different boxes. This actually arises in collecting baseball cards. Kolchin et al. [37] give some classical results.

3. *Coverage problems.* A sweeping generalization occurs in the world of coverage processes. Let $X$ be a compact metric space. Pick points $x_1, \ldots, x_k$ from a probability measure $\mu$ on $X$. For a fixed $\epsilon > 0$, what is the chance that the union of the $\epsilon$-balls about the $x_i$ cover $X$? Classical examples are dropping arcs on a circle (see Feller [26]) or caps on a sphere. Hall [31] and Janson [32] give many results. Aldous [1] has further literature. The widely studied k-SAT problem of computer science theory can be seen as a special case. Here, $X$ is the hyper-cube $\mathbb{Z}_2^d$, $k$ points are chosen randomly, and one wants to understand the chance that Hamming balls of a given radius cover the cube.

### 11.3. The birthday problem

Many natural variations have been studied: one may fix a graph, color the vertices with $k$ colors chosen with probabilities $p_i$ for $1 \leq i \leq k$, and ask how many edges have the same color. The same kind of question may be asked for hyper-graphs. Some of these problems are treated by Aldous [1], Arratia-Gordon-Goldstein [4], or Barbour-Holst-Janson [8]. We mention here a lesser-known development. Camarri and Pitman [12] have shown that limits other than the Poisson can arise if non-uniform probabilities are used. The limiting measures that arise there are natural (related to the Mittag-Leffler function) and well worth further development. It is also natural to study how these limiting regimes occur in the graph variants discussed above. Developing Stein's method for these cases seems like a natural project.

We close this brief section with a pointer to infinitely many other open problems. David Aldous has written a remarkable book called *Probability Approximations via the Poisson Clumping Heuristic*. This is really a collection of insights and conjectures, almost all related to Poisson approximation. Almost all of hundreds of problems are open.